\documentclass[a4paper,12pt, british, reqno]{amsart}

\usepackage{babel}
\usepackage{mathabx}
\usepackage{amssymb}
\usepackage{amsfonts}
\usepackage{enumitem}
\usepackage{hyperref}
\usepackage[utf8]{inputenc}
\usepackage{newunicodechar}
\usepackage{mathtools}
\usepackage{varioref}
\usepackage[centering]{geometry}
\usepackage[arrow,curve,matrix]{xy}

\usepackage{colortbl}
\usepackage{graphicx}
\usepackage{tikz}
\usepackage{tikz-cd}
\usepackage{capt-of}

\usepackage{mathrsfs}

\usepackage{multirow}
\usepackage{multicol}
\usepackage{longtable} 
\usepackage{caption}
\newcolumntype{M}[1]{>{\centering\arraybackslash}m{#1}} 

\definecolor{linkred}{rgb}{0.7,0.2,0.2}
\definecolor{linkblue}{rgb}{0,0.2,0.6}

\numberwithin{figure}{section}

\setdescription{labelindent=\parindent, leftmargin=2\parindent}
\setitemize[1]{labelindent=\parindent, leftmargin=2\parindent}
\setenumerate[1]{labelindent=0cm, leftmargin=*, widest=iiii}

\DeclareFontFamily{OMS}{rsfs}{\skewchar\font'60}
\DeclareFontShape{OMS}{rsfs}{m}{n}{<-5>rsfs5 <5-7>rsfs7 <7->rsfs10 }{}
\DeclareSymbolFont{rsfs}{OMS}{rsfs}{m}{n}
\DeclareSymbolFontAlphabet{\scr}{rsfs}
\DeclareSymbolFontAlphabet{\scr}{rsfs}

\DeclareMathOperator{\Ad}{Ad}
\DeclareMathOperator{\Aut}{Aut}

\DeclareMathOperator{\codim}{codim}

\DeclareMathOperator{\Pic}{Pic}

\DeclareMathOperator{\reg}{reg}

\DeclareMathOperator{\sing}{sing}
\DeclareMathOperator{\Spec}{Spec}
\DeclareMathOperator{\Sym}{Sym}

\DeclareMathOperator{\Gr}{Gr}

\DeclareMathOperator{\Exc}{Exc}

\DeclareMathOperator{\Cl}{Cl}
\DeclareMathOperator{\Cox}{Cox}

\DeclareMathOperator{\wt}{wt}

\newcommand{\sO}{\scr{O}}

\newcommand{\cB}{\mathcal B}
\newcommand{\cC}{\mathcal C}

\newcommand{\cH}{\mathcal H}

\newcommand{\cL}{\mathcal L}
\newcommand{\cM}{\mathcal M}

\newcommand{\cO}{\mathcal O}

\newcommand{\cU}{\mathcal U}

\newcommand{\cX}{\mathcal X}
\newcommand{\cY}{\mathcal Y}
\newcommand{\cZ}{\mathcal Z}
\newcommand{\0}{\mathcal O}
\newcommand{\bA}{\mathbb{A}}

\newcommand{\bC}{\mathbb{C}}

\newcommand{\bP}{\mathbb{P}}
\newcommand{\bQ}{\mathbb{Q}}
\newcommand{\bR}{\mathbb{R}}

\newcommand{\bT}{\mathbb{T}}

\newcommand{\bV}{\mathbb{V}}

\newcommand{\bZ}{\mathbb{Z}}

\newcommand{\bfa}{\mathbf{a}}

\newcommand{\bfc}{\mathbf{c}}

\newcommand{\bfm}{\mathbf{m}}
\newcommand{\bfM}{\mathbf{M}}
\newcommand{\bfN}{\mathbf{N}}
\newcommand{\bfP}{\mathbf{P}}

\newcommand{\bfu}{\mathbf{u}}
\newcommand{\bfv}{\mathbf{v}}

\newcommand{\fg}{\mathfrak{g}}
\newcommand{\fh}{\mathfrak{h}}

\newcommand{\fm}{\mathfrak{m}}

\newcommand{\fp}{\mathfrak{p}}

\newcommand{\ft}{\mathfrak{t}}
\newcommand{\fv}{\mathfrak{v}}

\newcommand{\aS}{{\sf S}}

\theoremstyle{plain}
\newtheorem{thm}{Theorem}[section]

\newtheorem{conjecture}[thm]{Conjecture}
\newtheorem{cor}[thm]{Corollary}
\newtheorem{defn}[thm]{Definition}

\newtheorem{lem}[thm]{Lemma}

\newtheorem{prop}[thm]{Proposition}

\theoremstyle{remark}

\newtheorem{c-n-d}[thm]{Claim and Definition}

\newtheorem{example}[thm]{Example}

\newtheorem{rem}[thm]{Remark}
\newtheorem{question}[thm]{Question}
\newtheorem*{rem-nonumber}{Remark}

\numberwithin{equation}{thm}

\setlist[enumerate]{label=(\thethm.\arabic*), before={\setcounter{enumi}{\value{equation}}}, after={\setcounter{equation}{\value{enumi}}}}

\newcommand{\factor}[2]{\left. \raise 2pt\hbox{$#1$} \right/\hskip -2pt\raise -2pt\hbox{$#2$}}

\author{Baohua Fu}
\address{Baohua Fu, State Key Laboratory of Mathematical Sciences, Morningside Center of Mathematics, Academy of Mathematics and Systems Science, Chinese Academy of Sciences, Beijing 100190, China;   and School of Mathematical Sciences, University of Chinese Academy of Sciences, Beijing, China}
\email{\href{bhfu@math.ac.cn}{bhfu@math.ac.cn}}
\urladdr{\href{http://www.math.ac.cn/people/fbh/}{http://www.math.ac.cn/people/fbh/}}

\author{Jie Liu} %
\address{Jie Liu, Institute of Mathematics, Academy of Mathematics and Systems Science, Chinese Academy of Sciences, Beijing, 100190, China}
\email{\href{jliu@amss.ac.cn}{jliu@amss.ac.cn}}
\urladdr{\href{http://www.jliumath.com}{http://www.jliumath.com}}


\makeatletter
\@namedef{subjclassname@2020}{2020 Mathematics Subject Classification}
\makeatother

\title[]{Symplectic singularities arising from algebras of symmetric tensors}
\date{\today}

\makeatletter

\hypersetup{
	pdfauthor={\authors},
	pdftitle={\@title},
	pdfsubject={\@subjclass},
	pdfkeywords={\@keywords},
	pdfstartview={Fit},
	pdfpagemode={UseOutlines},
	pdfpagelayout={OneColumn},
	bookmarks,
	colorlinks,
	linkcolor=linkblue,
	citecolor=linkred,
	urlcolor=linkred
}
\makeatother

\DeclareMathOperator{\Proj}{Proj}

\begin{document}

\begin{abstract}
	The algebra of symmetric tensors $S(X)\coloneqq H^0(X, \aS^{\bullet} T_X)$ of a projective manifold $X$ leads to a natural dominant  affinization morphism 
    \[
    \varphi_X\colon  T^*X \longrightarrow \cZ_X\coloneqq \Spec S(X).
    \]
    It is shown that $\varphi_X$ is birational if and only if $T_X$ is big. We prove that if $\varphi_X$ is birational, then $\cZ_X$ is a symplectic variety endowed with the Schouten--Nijenhuis bracket if and only if $\bP T_X$ is of Fano type, which is the case for smooth projective toric varieties, smooth horospherical varieties with small boundary, and the quintic del Pezzo threefold. These give examples of a distinguished class of conical symplectic varieties, which we call symplectic orbifold cones.
\end{abstract}

	\maketitle
	\tableofcontents
	
	\section{Introduction}

The cotangent bundle $T^*X$ of a projective manifold $X$ carries a natural holomorphic symplectic structure, making it a non-compact holomorphic symplectic manifold. Although the symplectic geometry of such spaces has been studied for specific examples \cite{Hwang2015}, the general case remains intriguing. In this paper, our aim is to provide a new perspective on the symplectic geometry of $T^*X$ by linking it to a novel construction of symplectic singularities.
    	
	\subsection{Motivation and question}
	
	A normal variety $Z$ is \emph{symplectic} in the sense of Beauville \cite{Beauville2000a} if there exists a symplectic form $\omega$ over the smooth locus $Z_{\reg}$ of $Z$ such that for any resolution $W\rightarrow Z$, the $2$-form $\omega$ extends to a holomorphic $2$-form on $W$. Moreover, a resolution $W\rightarrow Z$ of a symplectic variety $Z$ is called a \emph{symplectic resolution} if $\omega$ extends to a symplectic form on $W$ \cite{Fu2006}. 
	
	Let $X$ be a projective manifold, and let $\bP T_X$ be the projectivization in the sense of Grothendieck. The \emph{algebra of symmetric tensors} $S(X)$ of $X$ is the graded $\bC$-algebra defined as:
	\begin{equation}
		\label{eq.Sym-Alg}
		S(X) := \bigoplus_{p \geq 0} H^0(X, \aS^p T_X) = \bigoplus_{p \geq 0} H^0(\mathbb{P} T_X,\sO_{\bP {T_X}}(p)).
	\end{equation}
	This $S(X)$ is an integrally closed domain, though not necessarily finitely generated (over $\bC$) in general. Let $\mathcal{Z}_X \coloneqq  \operatorname{Spec} S(X)$ be the associated normal affine scheme.
	
	Note that $S(X)$ is isomorphic to the ring $\mathcal{O}(T^*X)$ of regular functions on $T^*X$ so that the symplectic structure on $T^*X$ induces a Poisson structure on $S(X)$ via the \emph{Schouten-Nijenhuis bracket} and we have a natural dominant affinization morphism
	\[
	\varphi_X\colon T^*X \longrightarrow \Spec \cO(T^*X) = \Spec S(X)\eqqcolon \mathcal{Z}_X.
	\]
	Despite its elementary construction, the pair $(\mathcal{Z}_X, \varphi_X)$ remains poorly understood in the general case (cf. \cite{beauville2024algebra}).
	
	\begin{defn}
		\label{d.symplectic-orbifold-cone}
		For a projective manifold $X$, we call $\mathcal{Z}_X$ a symplectic orbifold cone if $\varphi_X$ is birational and $\mathcal{Z}_X$ is a symplectic variety under the Schouten-Nijenhuis bracket.
	\end{defn}
	
	We refer the reader to Definition \ref{d.orbifold-cone} and Remark \ref{r.name-orbifold-cone} for terminology justification. In this paper, we are interested in the following questions.
	
	\begin{question}
		\label{q.Main-question}
		Let $X$ be a projective manifold with morphism $\varphi_X: T^*X \to \mathcal{Z}_X$.
		\begin{enumerate}
			\item\label{i.bir} When is $\varphi_X$ birational?
			\item\label{i.symp} When is $\mathcal{Z}_X$ a symplectic orbifold cone?
			\item\label{i.sympres} When is $\varphi_X$ a symplectic resolution?
		\end{enumerate}
	\end{question}

	\subsection{Main results}
	
	Recall that a normal projective variety $Y$ is \emph{of Fano type} if there exists an effective $\mathbb{Q}$-divisor $\Delta$ over $Y$ such that $-(K_Y + \Delta)$ is ample and $(Y, \Delta)$ is klt \cite[Definition 2.8]{Kollar2013}. For a vector bundle $E$ over a projective manifold, we say $E$ is \emph{big} (resp. \emph{nef}, \emph{semi-ample}) if $\sO_{\mathbb{P}E}(1)$ is. Our main result is the following.
	
	\begin{thm}
		\label{t.criterion-contact-symplectic}
		Let $X$ be a projective manifold.
		\begin{enumerate}
			\item\label{i.bira-big} The map $\varphi_X$ is birational if and only if $T_X$ is big.
			\item The scheme $\mathcal{Z}_X$ is a symplectic orbifold cone if and only if $\mathbb{P} T_X$ is of Fano type. In particular, if $\cZ_X$ is a symplectic orbifold cone, then $X$ is of Fano type.
		\end{enumerate}
	\end{thm}
    
    This resolves Questions \ref{i.bir} and \ref{i.symp}, leaving Question \ref{i.sympres} open. Proposition \ref{p.Q3-CPConj} shows its direct connection to the Campana--Peternell conjecture on Fano manifolds with nef tangent bundles. { Moreover, we remark that the statement \ref{i.bira-big} is more-or-less standard and it is known that $T_X$ big is also equivalent to $\varphi_X$ generically finite \cite[\S\,4]{beauville2024algebra}}.

	We now present some examples of symplectic orbifold cones. First note that the bigness of tangent bundle is highly restrictive for projective manifolds (cf. \cite{Hsiao2015,HoeringLiu2023,KimKimLee2025,ShaoZhong2025}). Known examples include smooth projective horospherical varieties \cite[\S\,3B2]{Liu2023a} and the quintic del Pezzo threefold \cite[Theorem 1.2]{HoeringLiuShao2022}.
	
	Recall that for a connected reductive group $G$, a normal $G$-variety is \emph{horospherical} if it admits an open $G$-orbit forming a torus bundle over a rational homogeneous space (\S\,\ref{ss.def-horospherical}).  Such varieties have \emph{small \textup{(resp. \emph{divisorial})} boundary} when the complement of the open orbit, the \emph{boundary}, has codimension at least two (resp. codimension one). Prototypical examples are projective toric varieties (with divisorial boundary) and rational homogeneous spaces (with small boundary).
	
	We will explicitly analyze $\mathcal{Z}_X$ for three classes:
	\begin{itemize}
		\item Smooth projective toric varieties
		\item Smooth projective horospherical $G$-varieties with small boundary
		\item The quintic del Pezzo threefold
	\end{itemize}
	Each case requires distinct techniques, discussed separately below.

\subsubsection{Smooth projective toric varieties}
 For a smooth projective toric variety $X$ defined by a fan $\Sigma$, we have the natural exact sequence
 $$
  0\longrightarrow \bfM \stackrel{B}{\longrightarrow} \bZ^{\Sigma(1)}=\bZ^N \stackrel{A}{\longrightarrow} \Pic(X)\longrightarrow 0,$$
  where $\Sigma(1)$ is the set of one-dimensional rays in $\Sigma$ and $\bfM$ is the lattice of characters of the torus. By using the integer matrix $A$, we can associate $X$ with a toric hyper-K\"ahler variety, denoted by 
  $\cY(A,0,0)$ (\S\,\ref{ss.toric-hyperkaher}), which is known to be a symplectic variety by \cite{Bellamy2023}. Based on an explicit computation of $S(X)$, it turns out that $\cZ_X$ coincides with $\cY(A,0,0)$, namely we have:
	\begin{thm}
 \label{t.toric-case}
		Let $X$ be a smooth projective toric variety { with dimension at least two}. Then $\cZ_X$ is isomorphic to the toric hyper-K\"ahler variety $\cY(A,0,0)$. 
	\end{thm}

{ The unique one-dimensional smooth projective toric variety is $\bP^1$, and $\cZ_{\bP^1}$ is well-known to be the quadratic cone in $\bC^3$.} As an immediate application of Theorem \ref{t.criterion-contact-symplectic} and Theorem \ref{t.toric-case}, we obtain the following seemly-new result, which generalizes \cite[Theorem 5.9]{HausenSues2010}.

 \begin{cor}
     The projectivized tangent bundle $\bP T_X$ of a smooth projective toric variety $X$ is  of Fano type. 
 \end{cor}

\subsubsection{Horospherical varieties with small boundary}

For a smooth projective horospherical $G$-variety $X$, we will relate the Cox ring of $\bP T_X$ to that of a weak Fano manifold and then a characterization of varieties of Fano type via Cox rings obtained in \cite{GongyoOkawaSannaiTakagi2015} allows us to obtain the following result.

   	\begin{thm}
    \label{t.Fanotype-MDS-horospherical}
	    Let $X$ be a smooth projective horospherical $G$-variety. Then $X$ is of Fano type and $\bP T_X$ is a Mori dream space. 
	    If moreover $X$ has small boundary, then $\bP T_X$ is of Fano type.
	\end{thm}

We refer the reader to \S\,\ref{ss.Cox-MDS} for the definition of Mori dream spaces. It is unclear for us whether the last statement still holds for \emph{all} smooth projective horospherical varieties. 

According to the classification result obtained by B.~Pasquier in \cite{Pasquier2009}, any smooth projective horospherical $G$-variety $X$ with Picard number one has small boundary and thus $\cZ_X$ is a symplectic variety. For such $X$, the identity component of its automorphism group $\widetilde{G}=\Aut^0(X)$ is worked out in {\em loc. cit.}. In particular, the variety $X$ carries a natural action by a torus $\bT\cong \bC^{\times}$, which is transitive and free along the torus fibers of the horospherical open orbit, such that there exists a homomorphism $\widebar{G}\coloneqq G\times \bT \rightarrow \widetilde{G}$ with finite kernel. We consider the associated moment maps (\S\,\ref{ss.moment-map}) for the actions of $\widebar{G}$ and $\widetilde{G}$, respectively, on $T^*X$:
\[
\begin{tikzcd}[column sep=large, row sep=large]
    T^*X  \arrow[r,"\Phi_X^{\widetilde{G}}"]  \arrow[dr,"\Phi_X^{\widebar{G}}" swap]
        & \cM_X^{\widetilde{G}} \arrow[r, phantom, "\subset"] \arrow[d]
            & \widetilde{\fg}^* \arrow[d] \\
        & \cM_X^{\widebar{G}} \arrow[r, phantom, "\subset"]
            & \widebar{\fg}^*.
\end{tikzcd}
\]
Our next result relates $\cZ_X$ to the moment map images $\cM_X^{\widetilde{G}}$ and $\cM_X^{\widebar{G}}$. Note that $\bT$ acts naturally on both $\cZ_X$ and $\cM_X^{\widetilde{G}}$. The moment maps $\Phi_X^{\widetilde{G}}$ and $\Phi_X^{\widebar{G}}$ factor through $\varphi_X\colon  T^*X\rightarrow \cZ_X$ via the natural maps $\aS^{\bullet} \widetilde{\fg}\rightarrow S(X)$ and $\aS^{\bullet}\widebar{\fg} \rightarrow S(X)$.

\begin{prop}
\label{p.PicardnubmeroneHorospherical}
    Let $X$ be a smooth projective horospherical $G$-variety with Picard number one, and let $\widetilde{G}$ be the identity component of its automorphism group. Then the induced map $\cZ_X\rightarrow \cM_X^{\widetilde{G}}$ is a generically finite morphism satisfying the following commutative diagram
    \[
    \begin{tikzcd}
        \cZ_X \arrow[r] \arrow[d]
            & \cM_X^{\widetilde{G}} \arrow[d] \\
        \widetilde{\cM}_X \arrow[r]
            & \cM_X^{\widebar{G}},
    \end{tikzcd}
    \]
    where $\widetilde{\cM}_X\cong \cZ_X/\!\!/\bT$, $\cM_{X}^{\widebar{G}}\cong \cM_{X}^{\widetilde{G}}/\!\!/ \bT$ and $\widetilde{\cM}_X\rightarrow \cM_X^{\widebar{G}}$ is a finite morphism.
\end{prop}

 \subsubsection{The quintic del Pezzo 3-fold}
 
 The quintic del Pezzo $3$-fold (denoted by $V_5$) is a smooth codimension $3$ linear section of $\mathrm{Gr}(2,5)$ under the Plücker embedding. Recent work by A.~H\"oring and T.~Peternell \cite{HoeringPeternell2024} explicitly constructed a weak Fano model of $\mathbb{P} T_{V_5}$, proving that $\mathbb{P} T_{V_5}$ is of Fano type. By Theorem \ref{t.criterion-contact-symplectic}, the scheme $\mathcal{Z}_{V_5}$ is therefore a 
 six-dimensional symplectic variety. Building on their construction, we show that
 \begin{itemize}
 	\item the singular locus $Z \subset \mathcal{Z}_{V_5}$ is an irreducible codimension two subset, and
 	\item the transverse singularities along general points of $Z$ are of type $A_2$, and
 	\item $\mathcal{Z}_{V_5}$ admits a unique $\mathbb{Q}$-factorial terminalization.
 \end{itemize}
 This unique terminalization yields a symplectic resolution of $\cZ_{V_5}$ outside the vertex  (Proposition \ref{p.QFT-ZV5}), which allows us to show that the symplectic variety $\mathcal{Z}_{V_5}$ itself admits no symplectic resolution (Proposition \ref{p.SympResV5}).

 \subsection{Organization of the paper}

In \S\,\ref{s.preliminaries}, we review foundations of orbifold cones, moment maps, Cox rings and varieties of Fano type. In \S\,\ref{s.proof-of-main-thm}, we discuss the relation between Question \ref{i.sympres} and the Campana--Peternell conjecture, and we prove Theorem \ref{t.criterion-contact-symplectic} in \S\,\ref{ss.Proof-MainThm}. In \S\,\ref{s.toric-var}, we study $\cZ_X$ for $X$ a smooth projective toric variety, proving Theorem \ref{t.toric-case} in \S\,\ref{ss.Proof-Toric}. In \S\,\ref{s.horospherical-var}, we study $\cZ_X$ for $X$ a smooth projective horospherical $G$-variety with small boundary, proving Theorem \ref{t.Fanotype-MDS-horospherical} in \S\,\ref{ss.Proof-Horo} and Proposition \ref{p.PicardnubmeroneHorospherical} in \S\,\ref{ss.Proof-Picardone}. In \S\,\ref{s.V5}, we study $\cZ_X$ for $X$ the quintic del Pezzo threefold via \cite{HoeringPeternell2024}.

\subsection{Conventions and notations}

Throughout this paper,  we work over the field $\bC$ of complex numbers. All varieties and manifolds are always assumed to be irreducible.  Let $\bT^n$ denote the $n$-dimensional torus $(\bC^{\times})^n$ whose Lie algebra will be denoted by $\ft_n$. When $n=1$, we will drop the index $n$. For a smooth variety $X$, $T_X$ denotes its tangent bundle viewed as a vector bundle over $X$ while $T^*X$ denotes the total space of its cotangent bundle. For a vector bundle $E$, we denote by $\bP E$ its projectivization in the sense of Grothendieck and by $\textbf{P} E$ its projectivization in the geometric sense, i.e., $\textbf{P}E=\bP E^*$. We will always denote by $\pi\colon T^*X\rightarrow X$ and $\widebar{\pi}\colon \bP T_X\rightarrow X$ the natural projections, respectively. We use $\sim$ (resp. $\sim_\bQ$) to denote (resp. $\bQ$-)rational equivalence.
 
\section{Preliminaries}
\label{s.preliminaries}

	\subsection{Orbifold cones}
	\label{ss.orbicone}
	
	An affine variety $Z=\Spec R$ is called a \emph{cone} if there exists a positive grading $R=\oplus_{i\geq 0} R_i$ such that $R_0=\bC$ and $R$ is a finitely generated domain over $R_0$. The grading determines a natural $\bT$-action on $Z$ and the vertex $0\in Z$ defined by the maximal ideal $\fm\coloneqq \oplus_{i\geq 1}R_i$ is the unique fixed point of this $\bT$-action.
	
	\begin{example}
		\label{e.affine-cone}
		Let $Y$ be a normal variety and let $L$ be a $\bQ$-Cartier Weil divisor on $Y$. The relative spectrum $\cL$ and the section ring $R(L)$ are defined as follows:
		\begin{center}
			
		\end{center}
		\[
		\cL\coloneqq \Spec_Y \left(\bigoplus_{p\geq 0}\sO_Y(pL)\right)\quad \textup{and}\quad R(L)\coloneqq \bigoplus_{p\geq 0} H^0(Y,\sO_Y(pL)).
		\]
        Then the natural projection $p\colon \cL \rightarrow Y$ is an affine morphism \cite[Proposition 1.3.2.3]{ArzhantsevDerenthalHausenLaface2015}. We also denote by $\cL^{\times}$ the associated Seifert $\bC^\times$-bundle \cite{Kollar2004}:
    \[
        \cL^{\times}\coloneqq \Spec_Y \left(\bigoplus_{p\in \bZ} \sO_Y(pL)\right).
    \]
    Then $\cL^{\times}$ is the open subset of $\cL$ obtained by removing the zero section $\textbf{0}_{\cL}$ of $\cL$ \cite[14]{Kollar2004}. 
    
    This $R(L)$ is an integrally closed domain and it is isomorphic to the ring of regular functions over $\cL$. Denote by $\cC_{L}$ the normal affine scheme $\Spec R(L)$ and by $\bP(\cC_L)$ the projective scheme $\Proj R(L)$. Then there exists a natural dominant affinization morphism 
		\begin{equation}
			\label{eq.Affinization-Cone}
			\varphi_L\colon \cL\longrightarrow \cC_{L}.
		\end{equation}
		Moreover, $\cC_L$ is a cone if $H^0(Y,\sO_Y)=\bC$ and $R(L)$ is a finitely generated $\bC$-algebra.
	\end{example}
	
	The following fact will be used in the proof of Theorem \ref{t.criterion-contact-symplectic}.
	
	 \begin{lem}
		\label{l.propreness}
		Let $Y$ be a normal projective variety and let $L$ be a $\bQ$-Cartier Weil divisor on $Y$. 
		\begin{enumerate}
			\item\label{i.proper-semi-ample} The map $\varphi_L$ is proper if and only if $L$ is semi-ample.
			
			\item\label{i.bir-big} The map $\varphi_L$ is birational if and only if $L$ is big.
		\end{enumerate}
	\end{lem}
	
	\begin{proof}
		For \ref{i.proper-semi-ample}, first we assume that $L$ is Cartier. Then $\cL\rightarrow Y$ is a $\bC$-bundle. If $\varphi_L$ is proper, it follows that $\varphi_L^{-1}(0)$ is complete, where $0\in \cC_L$ is the vertex, and thus $\varphi_L^{-1}(0)=\textbf{0}_{\cL} $ since $\varphi_L$ is $\bT$-equivariant. This means exactly that $L$ is semi-ample. 
        
        On the other hand, suppose that $L$ is semi-ample and consider the contraction $\widebar{\varphi}_L\colon Y\rightarrow \bP(\cC_L)$. Then there exists an ample Cartier divisor $\widebar{L}$ on $\bP(\cC_L)$ such that $\widebar{\varphi}_L^*\widebar{L}=L$ and $\widebar{\varphi}_L^*\colon R(\widebar{L})\rightarrow R(L)$ is an isomorphism of graded $\bC$-algebras. Moreover, as $\cL=Y\times_{\bP(\cC_L)}\widebar{\cL}$, where $\widebar{\cL}$ is the relative spectrum of $\widebar{L}$, the map $\varphi_L$ is proper if and only if $\varphi_{\widebar{L}}$ is proper. The latter is true since $\varphi_{\widebar{L}}\colon \widebar{\cL}\rightarrow \cC_{\widebar{L}}=\cC_L$ is actually a weighted blow-up of $0\in \cC_L$ (cf. \cite[\S\,3.1]{Kollar2013})
		
		For the general case, since $L$ is $\bQ$-Cartier, we can choose a positive integer $m$ such that $L'\coloneqq mL$ is Cartier. Consider the following commutative diagram
		\[
		\begin{tikzcd}[row sep=large,column sep=large]
			\cL \arrow[r,"\varphi_L"] \arrow[d,"q'" left]
			& \cC_{L} \arrow[d,"q"] \\
			\cL' \arrow[r,"\varphi_{L'}"]
			& \cC_{L'}
		\end{tikzcd}
		\]
		where $\cL'$ is the relative spectrum of $L'$. Note that $q'$ is finite, so $\varphi_L$ is proper if $\varphi_{L'}$ is proper. Conversely, if $\varphi_L$ is proper, then we have 
		\[
		(\varphi_{L'}\circ q')^{-1}(0')=(q\circ\varphi_L)^{-1}(0')= \varphi_L^{-1}(0) = \textbf{0}_{\cL},
		\]
		where $0'\in \cC_{L'}$ is the vertex. This implies $\varphi_{L'}^{-1}(0')=\mathbf{0}_{\cL'}$, so $L'$ is semi-ample and thus $\varphi_{L'}$ is proper. Since $L$ is semi-ample if and only if $L'$ is, the conclusion follows.
		
		For \ref{i.bir-big}, by \cite[Lemma 7.2]{Cutkosky2014}, $L$ is big if and only if $\dim R(X)=\dim(Y)+1$. In particular, as $\dim R(L)=\dim \cC_L$, it follows that $L$ is big if $\varphi_L$ is birational. Now we suppose that $L$ is big and consider the following commutative diagram
		\[
		\begin{tikzcd}[row sep=large, column sep=large]
			\cL \arrow[d,"p" left] \arrow[r,"\varphi_L"]
			&  \cC_L \arrow[d, dashed] \\
			Y   \arrow[r,dashed, "\widebar{\varphi}_L"]
			&  \bP(\cC_L)
		\end{tikzcd}
		\]
		Since $\widebar{\varphi}_L\colon Y\dashrightarrow \bP(\cC_L)$ is birational and $\varphi_L$ is $\bT$-equivariant, it remains to show that for a general point $y\in Y$, the restriction of $\varphi_L$ to the fiber 
        \[
        \bC\cong \cL_y\coloneqq p^{-1}(y)
        \]
        is generically injective. By \cite[Corollary 2.2.10]{Lazarsfeld2004}, there exists some $m\gg 1$ such that $H^0(Y,\sO_Y(mL))\not=0$ and $H^0(Y,\sO_X((m+1)L))\not=0$. Let $s_i\in H^0(Y,\sO_Y(iL))$ be a non-zero element for $i=m$ and $m+1$. As $y$ is general, we may also assume that $s_i|_{\cL_y}\not\equiv 0$. For two non-zero points $v$ and $v'\in \cL_y$, there exists a constant $\lambda\in \bT$ such that $v'=\lambda v$. In particular, we have $s_i(v')=\lambda^i s_i(v)$. As $(m,m+1)=1$, it follows that $s_i(v)=s_i(v')$, $i\in \{m,m+1\}$, if and only if $\lambda=1$, i.e., $v=v'$.
	\end{proof}
	
	Orbifold cones, which generalize cones over Fano manifolds, are among the most important examples of cones; see \cite[\S,3.1]{Zhuang2024a} and \cite[\S\,3.1]{Kollar2013} for details.
	
	\begin{defn}
		\label{d.orbifold-cone}
		Let $Y$ be a normal projective variety and let $L$ be an ample $\bQ$-Cartier Weil divisor on $Y$. The affine cone $\cC_L$ is called an orbifold cone if there exists some $r\in \bQ_{>0}$ such that $-K_Y\sim_{\bQ} r L$.
	\end{defn}
	
	\begin{lem} [\protect{\cite[42]{Kollar2004} and \cite[Lemma 3.3]{Zhuang2024a}}]
		\label{l.sing-orbifold-cone}
		An orbifold cone $\cC_L$ is klt if and only if $Y$ itself is klt.
	\end{lem}

		\subsection{Conical symplectic varieties}
	
	A symplectic variety $(Z,\omega)$ is called \emph{conical} if $Z=\Spec (R)$ is an affine cone and $\omega$ is homogeneous with respect to the natural $\bT$-action; that is, there exists a non-zero integer $l\in \bZ$ such that $\lambda^*\omega = \lambda^l \omega$ for any $\lambda\in \bT$. In this case we call $l$ the \emph{weight} of $(Z,\omega)$. Symplectic orbifold cones form a particularly interesting class of conical symplectic varieties with weight $1$. An interesting question is to characterize symplectic orbifold cones among all conical symplectic varieties.
	
	{  Now we consider a symplectic resolution $\mu\colon W \to Z$ of a conical symplectic variety  with weight $1$. The $\bT$-action on $Z$ lifts to $W$, which allows us to use the Bia\l ynicki-Birula decomposition \cite{B-B}.  With an adaption of arguments in  \cite[Lemmas 3.5-3.7]{Fu03}, we have the following:}
   
\begin{prop}
\label{p.conicalsymp}
        Assume that $(Z,\omega)$ is a conical symplectic variety with weight $1$ and  $\mu\colon W \to Z$ is a symplectic resolution with the lifted symplectic structure $\widetilde{\omega}$. Then there exists an irreducible component $X$ of $\mu^{-1}(0)$ contained in the $\bT$-fixed locus $W^{\bT}$, and an open subset $U \subset W$ such that $(U, \widetilde{\omega})  \simeq (T^*X, \omega_s)$ as conical symplectic varieties, where $\omega_s$ is the standard symplectic structure on the cotangent bundle.
    \end{prop}

\subsection{Moment maps}
\label{ss.moment-map}

Let $G$ be a linear algebraic group and denote by $\fg$ its Lie algebra. Let $X$ be a smooth $G$-variety. Given an arbitrary point $x\in X$, the tangent map of the orbit map $\mu_x\colon G\rightarrow Gx$, $g\mapsto gx$, defines a linear map 
\[
d\mu_x|_e\colon \fg\longrightarrow T_x (Gx) \subset T_x X.
\]
Taking dual yields a linear map $T^*_x X\rightarrow \fg^*$ and thus a morphism $\Phi_X^G\colon T^*X\rightarrow \fg^*$, which is called the \emph{moment map} (associated to the $G$-action). By $\cM_X^G\subset \fg^*$ we denote the closure of its image. Denote the image of the natural map $\aS^{\bullet}\fg \rightarrow S(X)$ by $R$. Then $\cM_X^G\cong \Spec R$ and we have a natural commutative diagram
\[
\begin{tikzcd}[row sep=large, column sep=large]
    T^*X \arrow[r,"\varphi_X"] \arrow[dr,"\Phi_X^G", swap]
        & \cZ_X \arrow[d] \\
        & \cM_X^G.
\end{tikzcd}
\]
The projectivization of $\Phi_X^{G}$ by the natural $\bT$-action induces a rational map $\bP T_X\dashrightarrow \bP(\cM_X^G)$, which is exactly the map defined by the linear system induced by the linear subspace $\fg \subset H^0(X,T_X)\cong H^0(\bP T_X,\sO_{\bP T_X}(1))$ (\cite[\S\,2C]{Liu2023a}).

\begin{example}
\label{e.momentmap}
    Let $X=\bC^n$ be the $n$-dimensional complex affine space. Then we can identify $T^*X$ to $\bC^n\times (\bC^{n})^*$ with coordinates $(x_1,\dots,x_n,y_1,\dots,y_n)$. 
    \begin{enumerate}
        \item Assume that the torus $\bT$ acts on $\bC^n$ with weight $(a_1,\dots,a_n)$; that is,
        \[
        t\cdot (x_1,\dots,x_n)\coloneqq (t^{a_1}x_1,\dots,t^{a_n} x_n).
        \]
        Then the moment map $\Phi_X^{\bT}$ can be written in coordinates as follows:
        \[
        \bC^n\times (\bC^n)^* \rightarrow \ft^*\cong \bC,\quad (x_1,\dots,x_n,y_1,\dots,y_n)\mapsto \sum_{i=1}^n a_i x_i y_i.
        \]

    \item Assume that the vector group $\bV\coloneqq \bC^n$ acts on $\bC^n$ by translation; that is,
    \[
    (v_1,\dots,v_n)\cdot (x_1,\dots,x_n) \coloneqq (v_1+x_1,\dots,v_n+x_n).
    \]
    Then the moment map $\Phi_X^{\bV}$ is nothing but the second projection, i.e.,
    \[
    \bC^n\times (\bC^n)^* \rightarrow \fv^*\cong \bC^n,\quad (x_1,\dots,x_n,y_1,\dots,y_n)\mapsto (y_1,\dots,y_n).
    \]
    \end{enumerate}
\end{example}

Here are some properties of the moment map:

\begin{prop}
Let $G$ be a linear algebraic group and let $X$ be a smooth projective $G$-variety.
\begin{enumerate}
    \item\label{i1.momentmap} The moment map   $\Phi_X^G$ is generically finite if and only if $G$ acts on $T^*X$ with an open orbit. In this case, the image $\cM_X^G$ is the closure of a co-adjoint orbit in $\fg^*$.

    \item\label{i2.momentmap} The moment map   $\Phi_X^G$ is proper if and only if $X$ is a projective $G$-homogeneous variety.

    \item\label{i3.momentmap} Assume that $G$ is reductive. Then $\Phi_X^G$ is generically finite if and only if $X$ is isomorphic to $G/P$ and the map $\Phi_X^G$ is the Springer map.
\end{enumerate}
\end{prop}

\begin{proof}
It is well-known (cf. \cite[\S\,8.2]{Timashev2011}) that for a point $\alpha \in T^*X$, the kernel ${\rm Ker} d_\alpha \Phi_X^G$  is the skew-orthocomplement of the tangent space of 
$G$-orbit $G\cdot \alpha$ and the image $d_\alpha\Phi_X^G)$ is the annihilator of the stabilizer $\fg_\alpha$.  It follows that $\Phi_X^G$ is generically finite if and only if $G$ acts on $T^*X$ with an open orbit. In this case, the image $\cM_X^G$ contains an open $G$-orbit, hence it is the closure of the co-adjoint orbit. This proves \ref{i1.momentmap}. 

For \ref{i2.momentmap}, assume that $\Phi_X^G$ is proper, then $(\Phi_X^G)^{-1}(0)$ is complete. By the $\bT$-equivariance, we have $X\cong\textbf{0}_X \subseteq (\Phi_X^G)^{-1}(0)$, where $\textbf{0}_X$ is the zero section of $T^*X\rightarrow X$, which we claim to be an equality. Assume $(x, \xi) \in T^*X$ is contained in $ (\Phi_X^G)^{-1}(0)$ for some $0 \neq \xi \in T_x^*X$, then $(x, \lambda \xi)$ is also contained in  $(\Phi_X^G)^{-1}(0)$ for all $\lambda \in \bT$. As $(\Phi_X^G)^{-1}(0)$ is complete, this implies that $\lim_{\lambda \to \infty} (x, \lambda \xi)$ lies in $T^*X$, which is absurd. This proves that $\textbf{0}_X = (\Phi_X^G)^{-1}(0)$, hence we have a morphism $\widebar{\Phi}\colon  \bP T_X \to \mathbb{P}\fg^*$ such that $\sO_{\bP T_X}(1) = \widebar{\Phi}^*(\sO_{\bP \fg^*}(1))$ is globally generated, hence $T_X$ is globally generated by $\fg\subset H^0(X,T_X)$, which implies that $X$ is $G$-homogeneous.

For \ref{i3.momentmap},  by \cite[Theorem 8.17]{Timashev2011}, the map $\Phi_X^G$ is generically finite if and only if the $G$-action on $X$ has complexity $0$ and rank $0$, which implies that $X$ is a rational homogeneous space by \cite[Proposition 10.1]{Timashev2011}.
\end{proof}

If the map $\Phi_X^G$ is generically finite, then $\cM_X^G$ is the closure of a co-adjoint orbit in $\fg^*$, which carries a symplectic structure on its open orbit. 
In general the variety $\cM_X^G$ is non-normal. It is then natural to pose the following:
\begin{question}
    Let $G$ be a linear algebraic group and let $X$ be a smooth projective $G$-variety.
 If the moment map $\Phi_X^G$ is generically finite, when is the normalization of $\cM_X^G$ a symplectic variety? 
\end{question}

\subsection{Varieties of Fano type}
	   Recall that a $\bQ$-divisor $\Delta$ on a normal projective variety is called \emph{big} if there exists an ample $\bQ$-divisor $A$ and an effective $\bQ$-divisor $N$ such that $\Delta\sim_{\bQ} A + N$. A birational map $g\colon X\dashrightarrow X'$ is called a \emph{birational contraction} if $g^{-1}$ does not contract any divisor.
	
	\begin{lem}[\protect{\cite[Lemma-Definition 2.6 and Lemma 2.8]{ProkhorovShokurov2009}}]
		\label{l.properties-Fano-type}
		Let $X$ be a normal projective variety.  
		\begin{enumerate}
			\item\label{i1.Fanotype} The variety $X$ is of Fano type if and only if there exists a big effective $\bQ$-divisor $\Delta$ such that $(X,\Delta)$ is klt and $K_X+\Delta \sim_{\bQ} 0$.
			
			\item\label{i2.Fanotype-birationalcontraction} If $X$ is of Fano type and $g\colon X\dashrightarrow X'$ is a birational contraction to a normal projective variety, then so is $X'$.
			
			\item\label{i3.Fanotype} If $X$ is of Fano type and $g: X \to X'$ is a projective contraction to a normal projective variety, then so is $X'$.
		\end{enumerate}
	\end{lem}
    \begin{proof}
        By \cite{ProkhorovShokurov2009}, it remains to prove \ref{i2.Fanotype-birationalcontraction}. Let $\Delta$ be an effective $\bQ$-divisor such that $(X,\Delta)$ is klt and $K_X+\Delta\sim_{\bQ} 0$. Set $\Delta'\coloneqq g_*\Delta$. Then $K_{X'}+\Delta'\sim_{\bQ} 0$. Let $\mu\colon W\rightarrow X$ be a log resolution of both $X$ and $g$. Let $\nu\colon W\rightarrow X'$ be the induced birational morphism. Then the Negativity Lemma implies 
        \[
        \nu^*(K_{X'}+\Delta') = \mu^*(K_X+\Delta).
        \]
        In particular, since $(X,\Delta)$ is klt and every $\mu$-exceptional divisor is $\nu$-exceptional, the pair $(X',\Delta')$ is klt and hence $X'$ is of Fano type.
    \end{proof}

    We also need the following characterization of varieties of Fano type via their anti-canonical models.
    
    \begin{thm}[\protect{\cite[Theorem 1.1]{CasciniGongyo2013}}] \label{t.CasciniGongyo}
    Let $X$ be a $\bQ$-Gorenstein normal projective variety. Then $X$ is of Fano type if and only if $-K_X$ is big, $R(-K_X)$ is finitely generated and $\Proj R(-K_X)$ has klt singularities.
    \end{thm}

\subsection{Cox rings and Mori dream spaces}
\label{ss.Cox-MDS}

Let $X$ be a $\bQ$-factorial normal projective variety. Choose a subgroup  $\Gamma$ of Weil divisors on $X$ such that $\Gamma_{\bQ}\rightarrow \Cl(X)_{\bQ}$ is an isomorphism. Then we define the Cox ring of $X$ to be
\[
\Cox(X,\Gamma)=\bigoplus_{D\in \Gamma} H^0(X,\sO_X(D))
\]
The definition depends on a choice of the group $\Gamma$ and if $\Cl(X)$ is free, then we will take $\Gamma$ such that $\Gamma\rightarrow \Cl(X)$ is an isomorphism and simply write $\Cox(X)$ for it. This is the case, for instance, if $X$ is smooth and $H^1(X,\sO_X)=0$. We refer the reader to \cite{ArzhantsevDerenthalHausenLaface2015} or \cite[\S\,2]{HausenSues2010} for more details about Cox rings in the general case.

Recall that a $\bQ$-factorial normal projective variety is called a \emph{Mori dream space} if its Cox ring is finitely generated \cite{HuKeel2000}. According to \cite[Corollary 1.3.2]{BirkarCasciniHaconMcKernan2010}, a $\bQ$-factorial normal projective variety of Fano type is a Mori dream space. Moreover, if $X$ is a Mori dream space and $D$ is a divisor on $X$, then there exists a minimal model program (MMP for short) for $D$; that is, a sequence of birational contractions
\[
X\coloneq X_0 \dashrightarrow X_1 \dashrightarrow \dots \dashrightarrow X_{m-1} \dashrightarrow X_m
\]
with $D_i$ the strict transform of $D$ on $X_i$ such that 
\begin{itemize}
    \item the $X_i$'s are $\bQ$-factorial normal projective varieties, and

    \item the birational map $X_i\dashrightarrow X_{i+1}$ is either a $D_i$-negative divisorial contraction or a $D_i$-negative flip, and

    \item either $D_m$ is nef, or $X_m$ has a $D_m$-negative elementary  contraction of fiber type.
\end{itemize}

We end this section with the following result, which characterizes varieties of Fano type via their Cox ring. {  We remark that the following statement does not depend on the choice of $\Gamma$ by \cite[Remark 2.9]{GongyoOkawaSannaiTakagi2015}.}

\begin{thm}[\protect{\cite[Theorem 1.1]{GongyoOkawaSannaiTakagi2015}}]
\label{t.GOST-Cox-Fano}
    Let $X$ be a $\bQ$-factorial normal projective variety. Then $X$ is of Fano type if and only if $\Cox(X,\Gamma)$ is finitely generated and $\Spec \Cox(X,\Gamma)$ has only klt singularities.
\end{thm}
	
\section{Symplectic orbifold cones}

\label{s.proof-of-main-thm}

This section is dedicated to the proof of Theorem \ref{t.criterion-contact-symplectic}. We will establish a more general result for quasi-contact pairs (Proposition \ref{p.criterion-contact-to-symplectic}), beginning with some basic definitions.

\subsection{Contact manifolds}
\label{ss.Contact-mfd}

Recall that a \emph{contact manifold} $(Y,L)$ is a pair consisting of a smooth variety $Y$, a Weil divisor $L$ on $Y$ and an exact sequence of vector bundles
\[
0\rightarrow F\rightarrow T_Y \xrightarrow{\theta}\sO_Y(L)\rightarrow 0
\]
such that the induced map by Lie bracket
\[
[\cdot,\cdot]\colon F\times F\rightarrow T_Y \xrightarrow{\theta} \sO_Y(L)
\]
is non-degenerate at every point of $Y$; or equivalently the twisted form
\begin{equation}
	\label{e.non-vanishing-contact}
	\theta\wedge(d\theta)^{n-1}\in H^0(Y,\sO_Y(K_Y+nL))
\end{equation}
is nonwhere vanishing, where $\dim (Y)=2n-1$. In particular, we have 
\begin{equation}
	\label{eq.Can-Div-Contact}
	-K_Y\sim nL. 
\end{equation}
We call $\sO_Y(L)$ the \emph{contact line bundle} and the $\sO_Y(L)$-valued $1$-form $\theta$, viewed as a global section of $\sO_Y(L)\otimes \Omega_Y^1$, is called the \emph{contact form}. Finally, we remark that the degeneracy locus of $\theta$ has pure codimension one and hence the non-degeneracy of $\theta$ can be verified over any open subset $U\subset Y$ with $\codim(Y\setminus U)\geq 2$.

As in Example \ref{e.affine-cone}, denote by $p\colon \cL\rightarrow Y$ the natural projection. Then $\cL^{\times}\rightarrow Y$ is a $\bC^{\times}$-bundle equipped with a natural symplectic form induced by the contact form $\theta$ \cite[\S\,2.2]{KebekusPeternellSommeseWisniewski2000}, which thereby induces a Poisson structure on $R(L)$. Moreover, if $H^0(Y,\sO_Y)=\bC$ and $\cC_L$ is an affine cone, then the affinization morphism $\varphi_{L}\colon \cL \longrightarrow \cC_{L}$ sends the zero section $\textbf{0}_{\cL}$ to the vertex $0\in \cC_L$. 

\begin{example}[\protect{\cite[\S\,2.6]{KebekusPeternellSommeseWisniewski2000}}]
	\label{e.cotangentbundle}
	Let $X$ be a projective manifold. Then the projectivized tangent bundle $(\bP T_X, L)$ is a contact manifold, where $L$ is a Weil divisor such that $\sO_{\bP T_X}(L)\cong \sO_{\bP T_X}(1)$. The contact form $\theta$ is given by the following composition:
	\[
	T_{\bP T_X} \longrightarrow \widebar{\pi}^*T_X \rightarrow \sO_{\bP T_X}(1).
	\]
	Moreover, $\cL$ is isomorphic to the blow-up of $T^*X$ along its zero section $\textbf{0}_X$ so that $\cL^{\times}\cong T^*X\setminus \textbf{0}_X$ is an isomorphism as symplectic varieties. By \eqref{eq.Sym-Alg}, we have $S(X)\cong R(L)$ so that $\cZ_X\cong \cC_L$.
\end{example}

 \begin{rem}
	By \cite[Theorem 1.1]{KebekusPeternellSommeseWisniewski2000} and \cite[Corollary 4]{Demailly2002}, a contact projective manifold $(Y,L)$ is either a Fano manifold with Picard number one or isomorphic to $(\bP T_X,L)$ in Example \ref{e.cotangentbundle} for some projective manifold $X$. 
\end{rem}

\subsection{Campana--Peternell conjecture}

We now relate Question \ref{i.sympres} to the Campana–Peternell conjecture. Although this connection is well-known to experts, we provide a complete proof due to the absence of explicit references. We begin with the following classical examples.

	\begin{example}[Springer map]
	\label{e.Springer-map}
	Let $G$ be a semi-simple linear algebraic group and let $P\subset G$ be a parabolic subgroup. We denote by $\fg$ and $\fp$ their Lie algebras. Then the cotangent bundle $T^*(G/P)$ of the rational homogeneous space $G/P$ can be naturally identified to the homogeneous bundle $G\times_P \fp^{\perp}$. Associating to a pair $(g,v)\in G\times \fp^{\perp}$ the element $\textup{Ad}(g)\cdot v$ of $\fg^*$ defines a generically finite proper map $\Phi\colon T^*(G/P)\rightarrow \fg^*$, called the \emph{Springer map}, whose image $\widebar{\cO}_P$ is the closure of a nilpotent orbit \cite{Richardson1974}. Let 
	\[
	T^*(G/P)\stackrel{\widetilde{\Phi}}{\longrightarrow} \widetilde{\cO}_P\longrightarrow \widebar{\cO}_P
	\]
	be the Stein factorization of $\Phi$. Then $\widetilde{\Phi}$ is a birational projective morphism, and hence $\widetilde{\cO}_P=\cZ_{G/P}$ and $\varphi_{G/P}=\widetilde{\Phi}$ so that$\widetilde{\Phi}$ is a symplectic resolution. In particular $T_{G/P}$ is big and semi-ample.
\end{example}

\begin{conjecture}[\protect{\cite[Conjecture 11.1]{CampanaPeternell1991}}]
	Let $X$ be a Fano manifold. Then $T_X$ is nef if and only if $X$ is isomorphic to a rational homogeneous space $G/P$.
\end{conjecture}

J.~Wang's recent result \cite{Wang2024} shows that for a projective manifold $X$, $T_X$ is big and semi-ample precisely when $X$ is Fano with nef tangent bundle. This directly connects Question \ref{i.sympres} to the Campana--Peternell conjecture.

\begin{prop}
	\label{p.Q3-CPConj}
	Let $X$ be a projective manifold. Then $T_X$ is big and semi-ample if and only if $\varphi_X\colon T^*X\rightarrow \cZ_X$ is a symplectic resolution.
\end{prop}

\begin{proof}
	As in Example \ref{e.cotangentbundle}, we consider the contact manifold $(\bP T_X,L)$, where $L$ is a Weil divisor corresponding to $\sO_{\bP T_X}(1)$. By \eqref{eq.Sym-Alg}, we have the following commutative diagram:
	\[
	\begin{tikzcd}[column sep=large, row sep=large]
		\cL \arrow[d,"\mu" left] \arrow[r,"\varphi_L"]
		& \cC_L\coloneqq \Spec R(L) \arrow[d,"\cong"] \\
		T^*X \arrow[r,"\varphi_X"]
		& \cZ_X\coloneqq \Spec S(X),
	\end{tikzcd}
	\]
	where $\mu$ is the blow-up along the zero section $\textbf{0}_X$ of $T^*X\rightarrow X$ so that $\textbf{0}_{\cL}=\mu^{-1}(\textbf{0}_X)$. The result follows from Lemma \ref{l.propreness} since $T_X$ is semi-ample (resp. big) if and only if $L$ is semi-ample (resp. big) by definition.
\end{proof}

\begin{cor}
	\label{c.symp-resolution}
	Assume that $(Z,\omega)$ is a conical symplectic variety with weight $1$ and  $\mu\colon W \to Z$ is a symplectic resolution with the lifted symplectic structure $\widetilde{\omega}$. If $\mu^{-1}(0)$ is irreducible, then $W = T^*X$ for a projective manifold $X$. If additionally $\dim Z \leq 10$, then $X$ is a rational homogeneous space and $\mu$ is a Springer map.
\end{cor}
\begin{proof}
	As $\mu^{-1}(0)$ is irreducible, we have $\mu^{-1}(0) = X$ for a projective manifold $X$ by Proposition \ref{p.conicalsymp}. {  As $\mu$ is $\bT$-equivariant, the $\bT$-fixed loci $W^{\bT}$ is contained in $\mu^{-1}(0)=X$,   hence $W^{\bT}=X$.  By the Bia\l ynicki-Birula decomposition \cite{B-B},} this implies that $W \cong T^*X$. It then follows that $Z \cong \cZ_X$ and $T_X$ is semi-ample and big. By \cite{Kan17}, the variety $X$ is a rational homogeneous space if $\dim X \leq 5$, which concludes the proof.
\end{proof}

\subsection{Proof of Theorem \ref{t.criterion-contact-symplectic}}

\label{ss.Proof-MainThm}

   We now prove Theorem \ref{t.criterion-contact-symplectic}, basing on the following criterion:

\begin{prop}[\protect{\cite[Theorem 6]{Namikawa2001}}]
	\label{p.criteria-extending-symplectic-form}
	A normal variety $Z$ is symplectic if and only if it has rational Gorenstein singularities and admits a symplectic form on some open $U\subset Z_{\reg}$ with $\codim(Z\setminus U)\geq 2$.
\end{prop}

To simplify the exposition, we introduce the following notion.

\begin{defn}
	A quasi-contact pair $(Y,L)$ consists of a normal projective variety $Y$ with a Weil divisor $L$ such that $(Y_{\reg},L|_{Y_{\reg}})$ is a contact manifold.
\end{defn}

\begin{lem}
	\label{l.reduc-quasi-contact}
	Let $(Y,L)$ be a quasi-contact pair with $\dim(Y)=2n-1$ such that $Y$ is of Fano type. Then there exists a quasi-contact pair $(Y',L')$ such that
	\begin{enumerate}
		\item there exists a birational contraction $Y\dashrightarrow Y'$ inducing an isomorphism $R(L)\cong R(L')$, and
		
		\item $Y'$ is a Fano variety wit klt singularities such that $-K_{Y'}\sim nL'$.
	\end{enumerate}
\end{lem}

\begin{proof}
	Let $\mu\colon \widetilde{Y}\rightarrow Y$ be a small $\bQ$-factorization of $Y$ \cite[Corollary 1.4.4]{BirkarCasciniHaconMcKernan2010}; that is, $\mu$ is a small birational contraction and $\widetilde{Y}$ is $\bQ$-factorial. Set $\widetilde{L}\coloneqq \mu^*L$. Then the natural pull-back 
	\[
	\mu^*\colon H^0(Y,\sO_Y(mL))\rightarrow H^0(\widetilde{Y},\sO_{\widetilde{Y}}(m\widetilde{L}))
	\]
	is an isomorphism for any $m\in \bZ$ and hence it yields an isomorphism $\mu^*\colon R(L)\rightarrow R(\widetilde{L})$ of graded $\bC$-algebras, which satisfies the following commutative diagram
	\[
	\begin{tikzcd}[column sep=large, row sep=large]
		\widetilde{\cL} \arrow[d,dashed] \arrow[r,"\varphi_{\widetilde{L}}"]
		& \cC_{\widetilde{L}} \arrow[d,"\cong"] \\
		\cL   \arrow[r,"\varphi_L"]
		& \cC_{L}.
	\end{tikzcd}
	\]
	Thus, after replacing $(Y,L)$ by $(\widetilde{Y},\widetilde{L})$, we may assume that $Y$ is $\bQ$-factorial and so $Y$ is a Mori dream space by \cite[Corollary 1.3.2]{BirkarCasciniHaconMcKernan2010}.
	
	Since $Y$ is of Fano type, $-K_Y$ is big. In particular, we can run a $-K_Y$-MMP $g\colon Y\dashrightarrow Y_1$ to yield a $\bQ$-factorial normal projective variety $Y_1$ with $-K_{Y_1}$ big and nef (\S\,\ref{ss.Cox-MDS}). Then Lemma \ref{l.properties-Fano-type} says that $Y_1$ is again of Fano type and hence the Basepoint Free Theorem implies that $-K_{Y_1}$ is actually big and semi-ample.
	
	Denote by $Y'=\Proj R(-K_{Y_1})$ and by $h\colon Y_1\rightarrow Y'$ the crepant contraction. Denote $g_*L$ by $L_1$ and $h_*L_1$ by $L'$. Then $L'$ is ample. By the Negativity Lemma, we have
	\[
	g^*h^*L'=g^*L_1=L-\Delta,
	\]
	where $\Delta$ is an effective $g$-exceptional divisor. For any positive integer $m$, the push-forward yields an inclusion
	\[
	\psi_m\colon H^0(Y,\sO_Y(mL)) \longrightarrow H^0(Y',\sO_{Y'}(mL')).
	\]
	We claim that $\psi_m$ is actually surjective. Indeed, let $f\in K(Y)=K(Y')$ be a rational function such that $f\in H^0(\sO_{Y'},\sO_{Y'}(mL'))$, i.e., $\textup{div}(f)+mL'\geq 0$. Then we have
	\[
	\textup{div}(f)+mL = \textup{div}(f) + mg^*h^*L'+ m\Delta \geq  g^*h^*(\textup{div}(f) + mL') \geq 0,
	\]
	which implies $f\in H^0(Y,\sO_Y(mL))$. Hence, the inclusion $R(L)\rightarrow R(L')$ is actually an isomorphism of graded $\bC$-algebras. 
	
	Since $Y_1$ is of Fano type and $\bQ$-factorial, it has only klt singularities, which implies that $Y'$ is klt as $h$ is crepant. Note that $Y'$ is  of Fano type by Lemma \ref{l.properties-Fano-type} and since $h\circ g$ is a birational contraction, $(Y',L')$ is again a quasi-contact pair. This yields $-K_{Y'}\sim nL'$ by \eqref{eq.Can-Div-Contact} and we are done.
\end{proof}

\begin{rem}
	The variety $Y'$ is actually the \emph{anti-canonical model} $\Proj R(-K_Y)$ of $Y$.
\end{rem}

For a quasi-contact pair $(Y,L)$, the natural restriction gives $R(L)\cong R(L|_{Y_{\reg}})$, inducing a Poisson structure on $\cC_L$ (\S\,\ref{ss.Contact-mfd}). Motivated by Question \ref{q.Main-question} and Example \ref{e.cotangentbundle}, we prove the following result, which implies Theorem \ref{t.criterion-contact-symplectic}.

\begin{prop}
	\label{p.criterion-contact-to-symplectic}
	Let $(Y,L)$ be a quasi-contact pair such that $L$ is big and $\bQ$-Cartier. Then $\cC_L$ is symplectic under the natural Poisson structure if and only if $Y$ is of Fano type.
\end{prop}

	\begin{proof}
        First we assume that $Y$ is of Fano type. By Lemma \ref{l.reduc-quasi-contact}, after replacing $(Y,L)$ by $(Y', L')$, we may assume that $L$ is ample, so  $\varphi_L\colon \cL\rightarrow \cC_L$ is a proper birational contraction by Lemma \ref{l.propreness} and $\cC_L$ is klt by Lemma \ref{l.sing-orbifold-cone}, it follows that the Poisson structure on $\cC_L$ defines a symplectic form $\omega$ over an open subset $U\subset \cC_{L}$ with $\codim(\cC_L\setminus U)\geq 2$. { Then $\omega^n$ generates $K_{U}$, where $\dim(Y)=2n-1$, which implies that $K_{\cC_L}\sim 0$ is Cartier} and hence $\cC_L$ has canonical singularities. Then we conclude by Proposition \ref{p.criteria-extending-symplectic-form}.

		Next we assume that $\cC_L$ is symplectic. Note that $\cC_L$ is actually the orbifold cone associated to $(\bP(\cC_L),L')$, where $L'=(\widebar{\varphi}_L)_*L$ and $\widebar{\varphi}_L\colon Y\dashrightarrow \bP(\cC_L)$ is the projectivization of $\varphi_L$. Then Lemma \ref{l.sing-orbifold-cone} implies that $\Proj R(-K_Y)\cong \bP(\cC_L)$ is klt and hence $Y$ is of Fano type by Theorem \ref{t.CasciniGongyo}.
	\end{proof}
	
	\begin{proof}[Proof of Theorem \ref{t.criterion-contact-symplectic}]
		Let $L$ be a Weil divisor such that $\sO_{\bP T_X}(L)\cong \sO_{\bP T_X}(1)$. Then there exists also an isomorphism of symplectic varieties $\cL^{\times}\rightarrow T^*X\setminus \textbf{0}_X$ (Example \ref{e.cotangentbundle}). Consider the following commutative diagram:
		\[
		\begin{tikzcd}[row sep=large, column sep=large]
			T^*X \arrow[r,hookleftarrow] \arrow[dr,"\varphi_X",swap]
			& T^*X\setminus \textbf{0}_X \arrow[r,"\cong"] \arrow[d]
			& \cL^{\times} \arrow[r,hookrightarrow] \arrow[d]
			& \cL \arrow[dl,"\varphi_{L}"] \\
			& \cZ_X \arrow[r,"\cong"]
			& \cC_{L}.
			& 
		\end{tikzcd}
		\]
  The first statement in Theorem \ref{t.criterion-contact-symplectic} follows from Lemma \ref{l.propreness}. The second statement follows by applying Proposition \ref{p.criterion-contact-to-symplectic} to $(\bP T_X,L)$ and Lemma \ref{l.properties-Fano-type} to $\bP T_X\rightarrow X$.
	\end{proof}
	
	\begin{rem}
		\label{r.name-orbifold-cone}
  By Lemma \ref{l.reduc-quasi-contact}, the conical symplectic variety $\cZ_X=\cC_{L}$ is indeed an orbifold cone in the sense of Definition \ref{d.orbifold-cone} . This justifies our Definition \ref{d.symplectic-orbifold-cone}.
	\end{rem}
 
	\section{Symplectic singularities from toric varieties}
\label{s.toric-var}

We investigate $\cZ_X$ for smooth projective toric varieties $X$ in this section. We refer the reader to \cite{CoxLittleSchenck2011} for general facts on toric varieties.

    \subsection{The Cox ring of varieties with torus action}

    \label{ss.Cox-torus-action}

 We briefly recall the results obtained in \cite{HausenSues2010} about the Cox ring of normal complete varieties with an algebraic torus action (see also \cite[\S\,4.4]{ArzhantsevDerenthalHausenLaface2015}). Let $X$ be a $\bQ$-factorial normal projective variety such that $\Cl(X)$ is finitely generated and free. Assume that $X$ carries an effective algebraic torus action $\bT^d\times X\rightarrow X$. For a point $x\in X$, denote by $\bT^d_x\subset \bT^d$ the isotropy group at $x$ and consider the non-empty $\bT^d$-invariant open subset 
\[
X^{\circ}\coloneqq \{x\in X\,|\,\dim \bT^d_x=0\}\subset X.
\]
{ By \cite[Corollary 2]{Sumihiro1974}, the variety $X^{\circ}$ is covered by $\bT^d$-stable affine open subsets $U_i$ with geometric quotient $U_i\rightarrow V_i\coloneqq U_i/\bT^d$. Gluing these affine varieties $V_i$ in a natural way yields a geometric quotient $q\colon X^{\circ}\rightarrow X^{\circ}/\bT^d$ with an irreducible normal but possibly non-separated orbit space $X^{\circ}/\bT^d$.} It follows from \cite[Proposition 3.5]{HausenSues2010} that there exists a \emph{separation} $f\colon X^{\circ}/\bT^d\dashrightarrow Y$ to a variety $Y$, i.e., a rational map which is a local isomorphism in codimension one. Clearly both $Y$ and $f$ are unique up to isomorphisms in codimension one. Denote by $D_1,\dots,D_s$ those $\bT^d$-invariant prime divisors which have a finite generic isotropy group of order $>1$, and denote by $E_1,\dots, E_m$ the $\bT^d$-invariant prime divisors supported in $X\setminus X^{\circ}$. Then after suitably shrinking, we may assume that there exist prime divisors $C_0,\dots,C_r$ on $Y$ such that each inverse image $f^{-1}(C_j)$ is a disjoint union of prime divisors $C_{ij}$, where $1\leq j\leq n_i$, the map $f$ is an isomorphism over $Y\setminus(C_0\cup \dots \cup C_r)$ and all the $D_j$'s occur among the $D_{ij}\coloneqq q^{-1}(C_{ij})$. Denote by $l_{ij}\in \bZ$ the order of the generic isotropy group of $D_{ij}$. For a prime divisor $D$, we denote by $1_{D}\in H^0(X,\sO_X(D))$ the canonical section associated to $D$.

\begin{thm}[\protect{\cite[Theorem 1.2]{HausenSues2010}}]
\label{t.HS-Coxring}
    There is a graded injection $\Cox(Y)\rightarrow \Cox(X)$ of Cox rings and the assignments $S_k\mapsto 1_{E_{k}}$ and $T_{ij}\mapsto 1_{D_{ij}}$ induce an isomorphism of $\Cl(X)$-graded rings
    \[
    \Cox(X) \cong \Cox(Y)[S_1,\dots,S_m,T_{ij};0\leq i\leq r,1\leq j\leq n_i]/\langle T_i^{l_i}-1_{C_i};0\leq i\leq r\rangle,
    \]
    where $T_i^{l_i}\coloneqq T_{i1}^{l_{i1}}\dots T_{in_i}^{l_{in_i}}$ and the $\Cl(X)$-grading on the right-hand side is defined by associating to $S_k$ the class of $E_k$ and to $T_{ij}$ the class of $D_{ij}$.
\end{thm}

Note that \emph{a priori} the separation $Y$ may be non-projective and $\Cl(Y)$ may be non-free. For the definition of Cox rings in this general case, we refer the reader to \cite[\S\,2]{HausenSues2010}. However, we will only use Theorem \ref{t.HS-Coxring} in the case where both $X$ and $Y$ are projective manifolds with free divisor class group.

\subsection{Toric hyper-K\"ahler varieties}

\label{ss.toric-hyperkaher}

We recall the construction of toric hyper-K\"ahler varieties \cite{Namikawa2023}. We identify the cotangent bundle $T^*\bC^N$ with the affine space $\bC^{2N}$ with coordinates $(z_1,\dots,z_N$, $w_1,\dots,w_N)$. The natural symplectic form $\omega$ on $\bC^{2N}$ can be written as
\[
\omega\coloneqq \sum_{1\leq i\leq N} dz_i\wedge dw_i.
\]
Let $A=[\textbf{a}_1,\dots,\textbf{a}_N]$ be an $(N-n)\times N$-integer matrix. Then $A$ determines a $\bT^{N-n}$-action on $\bC^N$ and hence a $\bT^{N-n}$-action on $T^*\bC^N=\bC^{2N}$. More precisely, write $\textbf{t}\coloneqq (t_1,\dots,t_{N-n})\in \bT^{N-n}=(\bC^{\times})^{N-n}$. Then the element $\textbf{t}$ acts on $\bC^{2N}$ as follows:
\[
\textbf{t}\cdot (z_1,\dots,z_N,w_1,\dots,w_N)\coloneqq (\textbf{t}^{\textbf{a}_1} z_1,\dots, \textbf{t}^{\textbf{a}_N} z_N, \textbf{t}^{-\textbf{a}_1} w_1,\dots, \textbf{t}^{-\textbf{a}_N} w_N),
\]
where $\textbf{t}^{\textbf{a}_i}\coloneqq t_1^{a_{1,i}}\cdots t_{N-n}^{a_{N-n,i}}$ with $\textbf{a}_i=(a_{j,i})_{1\leq j\leq N-n}$. Then the moment map of the $\bT^{N-n}$-action on $T^*\bC^N$ is given by
\[
\Phi\colon T^*\bC^N=\bC^{2N} \longrightarrow \ft_{N-n}^*=\bC^{N-n},\quad (z_1,\dots,z_N,w_1,\dots,w_N)\longmapsto \sum_{1\leq i\leq N}\textbf{a}_i z_i w_i
\]
From now on we assume that 
\begin{itemize}
    \item the map $\bZ^N\xrightarrow{A} \bZ^{N-n}$ is surjective, and 

    \item there exists an $N\times n$-integer matrix $B$ such that no row of $B$ is zero and that the following sequence is exact:
\[
0\longrightarrow \bZ^{n} \stackrel{B}{\longrightarrow} \bZ^N \stackrel{A}{\longrightarrow} \bZ^{N-n} \longrightarrow 0.
\]
\end{itemize}
Let $\theta$ be a rational character of $\bT^{N-n}$ and $\xi\in \ft_{N-n}^*$. Then we define
\[
\cY(A,\theta,\xi)\coloneqq \Phi^{-1}(\xi)^{\theta}/\!\!/\bT^{N-n},
\]
which is called a \emph{toric hyper-K\"ahler variety} (aka \emph{hypertoric variety} in the literature). Here $\Phi^{-1}(\xi)^{\theta} \subset \Phi^{-1}(\xi)$ denotes the open subset of $\theta$-semi-stable points {  in the sense of geometric invariant theory, which can be explicitly described as follows (see for example \cite[p.466]{Nag21}:}
$$
\Phi^{-1}(\xi)^{\theta} =\left\{(\textbf{z}, \textbf{w}) \in \Phi^{-1}(\xi) \,\middle\vert\, \theta \in \sum_{i: z_i \neq 0} \bQ_{\geq 0} \bf{a}_i + \sum_{j: w_j\neq 0} \bQ_{\geq 0}(-\bf{a}_j) \right\}
$$
According to \cite[Proposition 2.5]{Bellamy2023}, the variety $\cY(A,\theta,\xi)$ is a symplectic variety equipped with the symplectic form induced by $\omega$.

\subsection{Proof of Theorem \ref{t.toric-case}}
\label{ss.Proof-Toric}
{ This subsection is devoted to prove Theorem \ref{t.toric-case}. The proof is based on Theorem \ref{t.HS-Coxring}: We give an explicit description of $S(X)$ and then compare it with the coordinate ring of $\cY(A,0,0)$.}

Let $X$ be an $n$-dimensional smooth projective toric variety with $n\geq 2$, given by a fan $\Sigma$ in a lattice $\bfN$, and $\bT^n\cong X_o\subset X$ the acting torus. There are two types of rays $\rho\in \Sigma(1)$: those with $-\rho\not\in \Sigma(1)$ and those with $-\rho\in \Sigma(1)$. Denote by $\Sigma_1$ a subset of $\Sigma(1)$ containing all the rays of the first type and one representative for every pair of the second type. Let $\bfv_{\rho}\in \bfN$ be the primitive generator of $\rho$. We define
\[
\Sigma_1^{\perp}\coloneqq \left\{\textbf{u}\coloneqq (u_{\rho})\in \bC^{\Sigma_1}\,\middle\vert\,\sum_{\rho\in \Sigma_1}u_{\rho} \bfv_{\rho} = \bf{0} \right\}.
\]

\subsubsection{Explicit description of $\Cox(\bP T_X)$}

 We briefly recall the following result obtained in \cite{HausenSues2010}, which gives an explicit description of $\Cox(\bP T_X)$.

\begin{thm}[\protect{\cite[Theorem 5.9]{HausenSues2010}}]
\label{t.Cox-Proj-tangent-toric}
    Let $X$ be an $n$-dimensional smooth projective toric variety given by a fan $\Sigma$ {  with $n \geq 2$}. Then the Cox ring of $\bP T_X$ is given by
    \[
    \Cox(\bP T_X) = \left. \bC[S_{\rho},T_{\tau};\, \rho\in \Sigma(1),\tau\in \Sigma_1\,] \middle/ \left\langle \sum_{\rho \in \Sigma_1} u_{\rho} S^{\rho} T_{\rho};\,\textbf{u}=(u_{\rho})\in \Sigma_1^{\perp} \right\rangle \right.
    \]
    where
    \[
    S^{\rho}\coloneqq
    \begin{cases}
        S_{\rho} S_{-\rho},  &  -\rho\in \Sigma(1);\\
        S_{\rho},            & \textup{else}.
    \end{cases}
    \]
\end{thm}

To simplify the notation, we set $\cX \coloneqq  \bP T_X$. There exists an open subset $\cX^{\circ}$ of $\cX$ with $\codim(\cX\setminus \cX^{\circ})\geq 2$ such that $\bT^n$ acts on $\cX^{\circ}$ with trivial isotropy subgroup and $\bP(\ft_n)$ is a separation of $\cX^{\circ}/\bT^n$ such that the projectivized moment map 
\[
\widebar{\Phi}\colon \cX=\bP T_X \dashrightarrow \bP(\ft_n)
\]
coincides with the composition
\[
\cX \stackrel{q}{\dashrightarrow} \cX^{\circ}/\bT^n \stackrel{f}{\dashrightarrow} \bP(\ft_n).
\]
Note that $\bfv_{\rho}$ can be regarded as an element in $\bT^n=\bfN\otimes_{\bZ} \bT$ and therefore an element in $\ft_n$. Denote by $C_{\rho}\subset \bP(\ft_n)$ the hyperplane defined by the natural quotient $\ft_n\rightarrow \ft_n/\bC \bfv_{\rho}$. Then the Cox ring of the projective space $\bP(\ft_n)$ can be naturally identified to
\[
    \left.\bC[U_{\rho};\,\rho\in \Sigma_1] \middle/ \left\langle \sum_{\rho \in \Sigma_1} u_{\rho} U_{\rho};\,\textbf{u}=(u_{\rho})\in \Sigma_1^{\perp} \right\rangle\right.,
\]
where $U_{\rho}$ corresponds to the canonical section $1_{C_{\rho}}$. { Note that the restriction of $\widebar{\Phi}$ to the fibers of $\widebar{\pi}$ over $X_o$ is an isomorphism,} so the pull-back $\widebar{\Phi}^*C_{\rho}$ contains a unique $\widebar{\pi}$-horizontal prime divisor, denoted by $H_{\rho}$. Then we have
\[
\widebar{\Phi}^*C_{\rho}=
\begin{cases}
    H_{\rho}+\widebar{\pi}^*D_{\rho}+\widebar{\pi}^*D_{-\rho}, 
        & -\rho\in \Sigma(1);\\
    H_{\rho}+\widebar{\pi}^*D_{\rho}, & \textup{else};
\end{cases}
\]
where $D_{\rho}$ is the toric divisor on $X$ corresponding to $\rho$. Then Theorem \ref{t.Cox-Proj-tangent-toric} follows from Theorem \ref{t.HS-Coxring} by identifying $S_{\rho}$ (resp. $T_{\tau}$) to the canonical section $1_{\widebar{\pi}^*D_{\rho}}$ (resp. $1_{H_{\tau}}$). 

\subsubsection{Characterization of elements in $S(X)$}

Now we aim to characterize the elements of $\Cox (\bP T_X)$ contained in the subring $S(X)$. To this end, we shall interpret $\Cox(\bP T_X)$ as the following ring:
    \begin{align*}
        R & \coloneqq \left. \bC[S_{\rho},T_{\rho};\, \rho\in \Sigma(1)\,] \middle/ \left\langle \sum_{\rho\in \Sigma_1} u_{\rho} S^{\rho} T_{\rho}, (u_{\rho})\in \Sigma_1^{\perp}; T_{\rho}+T_{-\rho}, \rho\in \Sigma(1)\setminus \Sigma_1 \right\rangle \right.\\
        & \cong \Cox(\bP T_X).
    \end{align*}
Let $\bfM$ be the lattice of characters of $\bT^n$. Set $N=|\Sigma(1)|$. Then we have the following short exact sequence
    \begin{equation}
    \label{e.Pic-toric}
        0\longrightarrow \bfM \stackrel{B}{\longrightarrow} \bZ^{\Sigma(1)}=\bZ^N \stackrel{A}{\longrightarrow} \Pic(X)\longrightarrow 0.
    \end{equation}
    Denote by $\ker(A)\subset \bZ^{\Sigma(1)}$ the kernel of the map $A$. Consider the following surjective map
    \[
    \bC[S_{\rho},T_{\rho};\rho\in \Sigma(1)] \longrightarrow R.
    \]
    Given $I'=(i'_{\rho}), I=(i_{\rho})\in \bZ_{\geq 0}^{\Sigma(1)}$, we denote by $S^{I'} T^I$ the monomial
    \[
    \prod_{\rho\in \Sigma(1)} S_{\rho}^{i'_{\rho}} T_{\rho}^{i_{\rho}}.
    \]
    By abuse of notation, we shall also use $S^{I'} T^I$ to denote its image in $R$.

    \begin{lem}
\label{l.char-monomials}
    We have $S^{I'}T^I\in S(X)$ if and only if $I'-I-\prescript{c}{}{I}\in \ker(A)$, where $\prescript{c}{}{I}=(\prescript{c}{}{i_{\rho}})$ is defined as follows:
    \[
    \prescript{c}{}{i_{\rho}}\coloneqq 
    \begin{cases}
        i_{-\rho} & \textup{if}\quad -\rho\in \Sigma(1);\\
        0         & \textup{else}.
    \end{cases}
    \]
\end{lem}

\begin{proof}
    By Theorem \ref{t.Cox-Proj-tangent-toric} and { the paragraph after it}, the element $S^{I'}T^I$ can be viewed as the canonical section of $\sO_{\bP T_X}(H_I+\widebar{\pi}^*D_{I'})$ with
    \[
    H_I \coloneqq \sum_{\rho\in \Sigma(1)} i_{\rho} H_{\rho},\quad D_{I'}\coloneqq \sum_{\rho\in \Sigma(1)} i'_{\rho} D_{\rho}
    \]
    where $H_{\rho}=H_{-\rho}$ if $\rho\in \Sigma(1)\setminus \Sigma_1$. As $\sO_{\bP T_X}(1) \sim H_{\rho} + \widebar{\pi}^*(D_{\rho} + D_{-\rho})$ for any $\rho\in \Sigma(1)$, where $D_{-\rho} \coloneqq 0$ if $-\rho\not\in \Sigma(1)$, the monomial $S^{I'} T^{I}$ corresponds to a section of
    \[
    \Sym^m T_X\otimes \sO_X(D_{I'}-D_{I}-D_{\prescript{c}{}{I}}),
    \]
    where $m=\sum i_{\rho}$. So $S^{I'}T^{I}$ is contained in $S(X)$ if and only if $D_{I'}\sim D_I + D_{\prescript{c}{}{I}}$, which is equivalent to $I'-I-\prescript{c}{}{I}\in \ker(A)$ by \eqref{e.Pic-toric}.
\end{proof}

{ Lemma \ref{l.char-monomials} means that $S(X)$ is the ring of $\bT^{N-n}$-invariants of $\Cox(\bP T_X)$ for the action defined by the matrix $A$.}
Before giving the proof of Theorem \ref{t.toric-case}, we discuss a concrete example to illustrate the idea.

\begin{example}
\label{e.blowupPn}
    Let $X\rightarrow \bP^n$ be the blow-up of $\bP^n$ at the point $p\coloneqq [0:\dots:0:1]$. Then $X$ is toric and its fan has ray generators given by $\textbf{v}_i=\textbf{e}_i$ for $1\leq i\leq n$, $\bfv_{n+1}=-\textbf{e}_1-\dots-\textbf{e}_{n}$ and $\textbf{v}_{n+2}=-\textbf{v}_{n+1}$. We choose $\Sigma_1$ to be $\{\bfv_1,\dots,\bfv_{n+1}\}$. Then $\Sigma_1^{\perp}=\bC\bfu$, where $\bfu=(1,\dots,1)$, and Theorem \ref{t.Cox-Proj-tangent-toric} implies
        \[
        \Cox(\bP T_X) \cong \left. \bC[x_1,\dots,x_{n+2},y_1,\dots,y_{n+1}] \middle/ \left\langle \sum_{i=1}^{n} x_i y_i + x_{n+1}x_{n+2} y_{n+1} \right\rangle \right..
        \]
        Now we add another variable $y_{n+2}$ and write $\Cox(\bP T_X)$ as
        \[
        R\coloneqq \left. \bC[x_1,\dots,x_{n+2},y_1,\dots,y_{n+2}] \middle/ \left\langle \sum_{i=1}^{n} x_i y_i + x_{n+1} x_{n+2} y_{n+1}, y_{n+1} + y_{n+2}\right\rangle \right..
        \]
        Note that the matrix $A$ associated to $X$ is 
        \[
        \begin{pmatrix}
            1 & \dots & 1 & 1 & 0 \\
            -1 & \dots & -1 & 0 & 1
        \end{pmatrix}.
        \]
        So the toric hyper-K\"ahler variety $\cY(A,0,0)$ associated to $A$ is the spectrum of the $\bT^2$-invariant subring of the following ring $R'$:
        \[
        R'\coloneqq \left.\bC[x'_1,\dots,x'_{n+2},y'_1,\dots,y'_{n+2}] \middle/  \left\langle \sum_{i=1}^{n+1} x'_i y'_i, \  -\sum_{i=1}^{n} x'_i y'_i + x'_{n+2} y'_{n+2}\right\rangle \right..
        \]
        We also note that $R'$ can be naturally identified to the following ring 
        \[
        \widetilde{R}'\coloneqq \left.\bC[x'_1,\dots,x'_{n+2},y'_1,\dots,y'_{n+2}] \middle/  \left\langle \sum_{i=1}^{n+1} x'_i y'_i, \  x'_{n+1} y'_{n+1} + x'_{n+2} y'_{n+2}\right\rangle \right..
        \]
        Let $\widetilde{R}$ be the ring defined as follows:
        \[
        \left. \bC[x_1,\dots,x_{n+2},y_1,\dots,y_{n+2}] \middle/ \left\langle \sum_{i=1}^{n} x_i y_i + x_{n+1} x_{n+2} y_{n+1}, \ x_{n+1} x_{n+2} (y_{n+1} + y_{n+2})\right\rangle \right..
        \]
        Then we get a natural commutative diagram
        \[
        \begin{tikzcd}[column sep=large, row sep=large]
            \widetilde{R}'  \arrow[r,"\widetilde{\phi}"] \arrow[d,"\cong" left] 
                & \widetilde{R} \arrow[d] \\
            R' \arrow[r,"\widebar{\phi}"]
                & R,
        \end{tikzcd}
        \]
        where $\widetilde{\phi}$ is defined as follows:
        \[
        \begin{cases}
            x'_i\mapsto x_i & 1\leq i\leq n+2; \\
            y'_i\mapsto y_i & 1\leq i\leq n; \\
            y'_{n+1}\mapsto x_{n+2} y_{n+1} & \\
            y'_{n+2} \mapsto x_{n+1} y_{n+2}. &
        \end{cases}
        \]
        Then $\widebar{\phi}^*\colon \Spec (R)\rightarrow \Spec(R')$ is birational and the $\bT^2$-action on $R'$ can be naturally extended to $R$ with
        \[
        \wt(y_{n+1}) = (-1,-1)\quad \textup{and} \quad \wt(y_{n+2}) = (-1,-1).
        \]
        Let $x^I y^J\coloneqq x_1^{i_1}\dots x_{n+2}^{i_{n+2}} y_1^{j_1}\dots y_{n+2}^{j_{n+2}}$ be a monomial. By Lemma \ref{l.char-monomials}, one can derive that $x^I y^J\in S(X)$ if and only if the following equalities hold:
        \[
        \begin{dcases}
            \sum_{k=1}^{n+1} (i_k-j_k) - j_{n+2} =0 \\
            \sum_{k=1}^n (j_k-i_k) + i_{n+2} - j_{n+1} - j_{n+2} = 0.
        \end{dcases}
        \]
        This yields
        \[
        i_{n+1} + i_{n+2} = 2(j_{n+1} + j_{n+2}).
        \]
        As $y_{n+2}=-y_{n+1}$ in $R$, after suitable replacement, we may assume that $i_{n+1}\geq j_{n+2}$ and $i_{n+2}\geq j_{n+1}$; that is, $x^I y^J$ is contained in the image $\widebar{\phi}$. This implies 
        \[
        S(X) \cong R^{\bT^2} = (R')^{\bT^2}.
        \]
\end{example}

\subsubsection{Auxiliary lemmas}

    Fix a basis of $\Pic (X)$ such that $\Pic(X)\cong \bZ^{N-n}$. Then we can write $A=[\textbf{a}_{\rho};\rho\in \Sigma(1)]$ and $\prescript{t}{}{A}=[\widebar{\bfa}_i;1\leq i\leq N-n]$. Note that the linear map $B$ is given as 
    \[
    \bfM \longrightarrow \bZ^{\Sigma(1)},\quad \bfm\longmapsto (\langle \bfv_{\rho},\bfm\rangle)_{\rho\in \Sigma(1)},
    \]
    so the exactness of \eqref{e.Pic-toric} implies that the vector space
    \[
    \Sigma(1)^{\perp}\coloneqq \left\{\bfc=(c_{\rho})\in  \bC^{\Sigma(1)} \,\middle\vert\,\sum_{\rho\in \Sigma(1)} c_{\rho} \bfv_{\rho}=0\right\}
    \]
    is generated by the $\widebar{\bfa}_i$'s $(1\leq i\leq N-n)$. Moreover, for any $\rho\in \Sigma(1)\setminus \Sigma_1$, we define an element 
    \[
    \prescript{\rho}{}{\bfc}\coloneqq (\prescript{\rho}{}{c}_{\rho'})\in \bC^{\Sigma(1)} 
    \]
    with 
    \[
    \prescript{\rho}{}{c}_{\rho'} \coloneqq 
    \begin{cases}
        1 & \text{if}\quad \rho'=\rho\quad\text{or}\quad -\rho;\\
        0 & \text{else}.
    \end{cases}
    \]
    Then the exactness of \eqref{e.Pic-toric} implies that $\prescript{\rho}{}{\bfc}$ is contained in $\Sigma(1)^{\perp}$.

    \begin{lem}
    \label{l.Signa_1perpSigma(1)perp}
        The space $\Sigma(1)^{\perp}$ is generated by $\Sigma_1^{\perp}$ and the $\prescript{\rho}{}{\bfc}$'s, where $\Sigma_{1}^{\perp}\hookrightarrow \Sigma(1)^{\perp}$ is given as
        \[
        \bfu \longmapsto (\bfu,0,\dots,0).
        \]
    \end{lem}

    \begin{proof}
     This follows from the fact that $\bfv_{\rho}+\bfv_{-\rho}=0$ if $\rho\in \Sigma(1)\setminus \Sigma_1$.
    \end{proof}

    Now we introduce a ring $R'$ as follows:
    \[
    R'\coloneqq \left. \bC[S_{\rho},T^{\rho};\, \rho\in \Sigma(1)\,] \middle/ \left\langle \sum_{\rho\in \Sigma(1)} \mathbf{a}_{\rho} S_{\rho} T^{\rho} \right\rangle \right..
    \]
    Then we have $\Spec R'\cong \Phi^{-1}(0)$ and hence $R'$ is a domain by \cite[Lemma 4.7]{BellamyKuwabara2012}. On the other hand, consider the following natural morphism
    \[
    \phi \colon \bC[S_{\rho},T^{\rho};\rho\in \Sigma(1)] \longrightarrow \bC[S_{\rho},T_{\rho};\rho\in \Sigma(1)]
    \]
    with
    \[
        T^{\rho} \longmapsto 
        \begin{cases}
            S_{-\rho} T_{\rho}  & \textup{if}\quad -\rho\in \Sigma(1);\\
            T_{\rho}            & \textup{else}.
        \end{cases}
    \]
    
    \begin{lem}
    \label{l.birationalmodification}
        The morphism $\phi$ induces a morphism $\widebar{\phi}\colon  R' \rightarrow R$ such that the induced morphism $\widebar{\phi}^*\colon \Spec R \rightarrow \Spec R'$ is birational.
    \end{lem}

    \begin{proof}
        Firstly we show that $\widebar{\phi}$ is well-defined. To see this, consider the ring $\widetilde{R}$ defined as follows:
        \[
        \left. \bC[S_{\rho},T_{\rho};\, \rho\in \Sigma(1)\,] \middle/ \left\langle \sum_{\rho\in \Sigma_1} u_{\rho} S^{\rho} T_{\rho}, (u_{\rho})\in \Sigma_1^{\perp}; S^{\rho} (T_{\rho} + T_{-\rho}), \rho\in \Sigma(1)\setminus \Sigma_1 \right\rangle \right..
        \]
        Then we have a natural quotient homomorphism $\widetilde{R}\rightarrow R$. On the other hand, consider the ring $\widetilde{R}'$ defined as follows:
        \[
        \left. \bC[S_{\rho},T^{\rho};\, \rho\in \Sigma(1)\,] \middle/ \left\langle \sum_{\rho\in \Sigma_1} u_{\rho} S_{\rho} T^{\rho}, (u_{\rho})\in \Sigma_1^{\perp}; S_{\rho} T^{\rho} + S_{-\rho}T^{-\rho}, \rho\in \Sigma(1)\setminus \Sigma_1 \right\rangle \right..
        \]
        Then Lemma \ref{l.Signa_1perpSigma(1)perp} implies that $\widetilde{R}'$ is isomorphic to the following ring:
        \[
        \left. \bC[S_{\rho},T^{\rho};\, \rho\in \Sigma(1)\,] \middle/ \left\langle \sum_{\rho\in \Sigma(1)} c_{\rho} S_{\rho} T^{\rho};\, \bfc=(c_{\rho})\in \Sigma(1)^{\perp} \right\rangle \right..
        \]
        In particular, since $\Sigma(1)^{\perp}$ is generated by the $\widebar{\bfa}_i$'s $(1\leq i\leq N-n)$, the ring $\widetilde{R}'$ is canonically identified to $R'$ and we have the following commutative diagram
        \[
        \begin{tikzcd}[row sep=large, column sep=large]
            \widetilde{R}' \arrow[r,"\widetilde{\phi}"] \arrow[d,"\cong" left]
                & \widetilde{R} \arrow[d] \\
            R' \arrow[r,"\widebar{\phi}"]
                & R,
        \end{tikzcd}
        \]
        where $\widetilde{\phi}$ is the natural morphism induced by $\phi$. Hence $\widebar{\phi}$ is well-defined.

        Note that $S_{\rho}\not=0$ in $R$ and and hence so is in $R'$ for any $\rho\in \Sigma(1)$. In particular, since $R$ and $R'$ are domains, the rational map $(\widebar{\phi}^*)^{-1}\colon \Spec R' \dashrightarrow \Spec R$ defined by sending $T_{\rho}$ to $T^{\rho}/S_{-\rho}$ if $\rho\in \Sigma(1)\setminus \Sigma_1$ shows that $\widebar{\phi}^*$ is indeed birational.
    \end{proof}

\subsubsection{Proof of Theorem \ref{t.toric-case}}

The torus $\bT^{N-n}$ acts on $R'$ with $\wt(S_{\rho})=\bfa_{\rho}$ and $\wt(T^{\rho})=-\bfa_{\rho}$ for any $\rho\in \Sigma(1)$. Note that $R'$ is a subalgebra of $R$ by Lemma \ref{l.birationalmodification}. Then the $\bT^{N-n}$-action can be extended to $R$ such that $\wt(T_{\rho})=-\bfa_{\rho}-\bfa_{-\rho}$ for any $\rho\in \Sigma(1)\setminus \Sigma_1$. An easy computation shows that a monomial $S^{I'} T^I\in R$, where $I'=(i'_{\rho}), I=(i_{\rho})\in \bZ_{\geq 0}^{\Sigma(1)}$,  has weight $\textbf{0}$ if and only if $I'-I - \prescript{c}{}{I}\in \ker(A)$ and the latter is equivalent to $S^{I'} T^I\in S(X)$ by Lemma \ref{l.char-monomials}. So it is enough to show that every $\bT^{N-n}$-invariant monomial $S^{I'} T^{I}$ in $R$ is contained in $R'$. 

Let $S^{I'}T^{I}$ be a $\bT^{N-n}$-invariant element in $R$. Then $I'-I-\prescript{c}{}{I}\in \ker(A)$, or equivalently $I'-I-\prescript{c}{}{I}\in \textup{Im}(B)$. In particular, as $\bfv_{\rho}=-\bfv_{-\rho}$ for any $\rho\in \Sigma(1)\setminus \Sigma_1$, we get 
\[
i'_{\rho} - i_{\rho} - i_{-\rho} = -i'_{-\rho} + i_{\rho} + i_{-\rho}
\]
for any $\rho$ such that $-\rho\in \Sigma(1)$. In particular, as $T_{\rho}=-T_{-\rho}$ in $R$, we may assume that $i'_{\rho}\geq i_{-\rho}$ if $-\rho\in \Sigma(1)$ and hence the monomial $S^{I'} T^I$ is contained in the subring $R'$ of $R$ generated by the $S_{\rho}$'s and the $S_{-\rho} T_{\rho}=T^{\rho}$'s. \qed

\begin{rem}
    It is noteworthy that extending Theorem \ref{t.toric-case} to any toric variety $X$ with $H^0(X,\sO_X^\times)=\bC^{\times}$ and $\Pic(X)$ free poses no essential difficulties (see \cite[Theorem 4.4.1.3]{ArzhantsevDerenthalHausenLaface2015}). However, given that the primary focus of this article is to present a distinctive perspective for studying $\cZ_X$ or equivalently $S(X)$, we choose to concentrate on smooth projective varieties.
\end{rem}

\subsection{Example: blow-ups of projective spaces along points}

Note that $\cZ_{\bP^n}$ is the minimal nilpotent orbit closure $\overline{\0}_\textup{min} \subset \mathfrak{sl}_{n+1}$ consisting of traceless matrices of rank at most 1. The map $\Phi_{\bP^n}: T^*\bP^n \to \overline{\0}_\textup{min}$ is the Springer resolution, which collapses the zero section to the origin. This map can be explicitly described as follows: Let $\bT$ act on $\bC^{n+1} \times (\bC^{n+1})^*$ by 
$$
\lambda \cdot (x_1, \cdots, x_{n+1}, y_1, \cdots, y_{n+1}) = (\lambda x_1, \cdots, \lambda  x_{n+1}, \lambda^{-1} y_1, \cdots, \lambda^{-1} y_{n+1}).
$$
This $\bT$-action preserves the hypersurface $\mathcal{H}\coloneqq \{(\textbf{x}, \textbf{y}) | \sum x_iy_i =0\}$. Moreover, since $\bP^n$ is a toric variety given by the fan $\Sigma_{\bP^n}$, which has rays generated by $\bfv_i=\textbf{e}_i$, $1\leq i\leq n$, and $\bfv_{n+1}=-\textbf{e}_1-\dots-\textbf{e}_n$, then Theorem \ref{t.Cox-Proj-tangent-toric} implies that
$\Spec S(\bP^n) \cong \mathcal{H} /\!\!/ \bT$. The Springer resolution $\Phi_{\bP^n}$ is then given by
$$
({\bf x}, {\bf y}) \longmapsto (x_iy_j)_{1\leq i, j \leq n+1}.
$$

Now let $X$ be the blow-up of $\bP^n$ at the toric point $[0:\dots:0:1]$. We will describe in details the symplectic singularity $\cZ_X$ (see Example \ref{e.blowupPn}). The toric structure of $X$ gives a $\bT^2$-action on $\bC^{n+2}$ via
$$
(t_1, t_2) \cdot (x_1,\cdots, x_{n+2}) = (t_1t_2^{-1}x_1, \cdots, t_1t_2^{-1}x_n, t_1 x_{n+1}, t_2 x_{n+2}).
$$
The associated moment map is then given by
$$
\Phi\colon \bC^{2n+4} \to \bC^2, \quad (x_i, y_j)  \mapsto \left(\sum_{i=1}^n x_iy_i+x_{n+1}y_{n+1}, -\sum_{i=1}^n x_iy_i+x_{n+2}y_{n+2}\right).
$$
The coordinate ring of the quotient $\Phi^{-1}(0)/\!\!/\bT^2$ is generated by the following monomials:
$$
x_iy_j\, (1 \leq i, j \leq n),\quad  x_iy_{n+1}x_{n+2}\, (1\leq i\leq n),\quad y_ix_{n+1}y_{n+2}\, (1\leq i\leq n).
$$
It follows from Theorem \ref{t.toric-case} that $\cZ_X \cong \Phi^{-1}(0)/\!\!/\bT^2 =\cY(A,0,0)$, which is isomorphic to the following
\[
\left\{ (z_{ij})_{1\leq i,j\leq n+1}\,\left\vert\, \sum_{i=1}^{n+1} z_{ii} =0, 
\quad {\rm rk} \begin{pmatrix} 
z_{1,1} & \cdots & z_{1, n} & z_{1,n+1} \\
\vdots & \vdots  & \vdots & \vdots \\
z_{n,1} & \cdots & z_{n, n} & z_{n,n+1} \\
z_{n+1, 1}& \cdots & z_{n+1, n} & -z_{n+1, n+1}^2
\end{pmatrix}  \leq 1 \right\}\right..
\]
Note that when $n=2$, this example is studied in \cite[Theorem 4.8]{Nag21}. The singular locus of $\cZ_X$ is given by 
$$
S\coloneqq \left.\left\{ C= \begin{pmatrix} 
z_{1,1} & \cdots & z_{1,n} & 0 \\
\vdots & \vdots  & \vdots  & \vdots \\
z_{n, 1}& \cdots & z_{n, n} & 0 \\
0 & \cdots & 0 & 0
\end{pmatrix}\, \right\vert\,  {\rm tr}(C)=0, {\rm rk}(C) \leq 1 \right\}
$$

Note that $S$ is isomorphic to the minimal nilpotent orbit closure in $\mathfrak{sl}_n$. Write
$S = S_0 \cup \{0\}$ as the union of rank 1 (resp. 0) traceless matrices. Note that $S$ has codimension two in $\cZ_X$ and the transverse singularity of $\cZ_X$ along any point of $S_0$ is an $A_1$-singularity. The blow-up $X \to \bP^n$ induces a natural morphism
\[
H^0(T^*X, \sO_{T^*X}) \longrightarrow H^0(T^* \bP^n, \sO_{T^*\bP^n}),
\]
which gives a map $\nu\colon   \overline{\0}_{\textup{min}} \to \cZ_X $.  In terms of matrices, this map is given by
$$
z\coloneqq \begin{pmatrix} 
z_{1,1} & \cdots & z_{1,n+1} \\
\vdots & \vdots  & \vdots \\
z_{n+1, 1}& \cdots & z_{n+1, n+1}
\end{pmatrix}  \mapsto \nu(z)\coloneqq 
\begin{pmatrix} 
z_{1,1} & \cdots & z_{1,n} & -z_{n+1, n+1} z_{1,n+1} \\
\vdots & \vdots  & \vdots &  \vdots \\
z_{n,1} & \cdots & z_{n, n} & -z_{n+1, n+1} z_{n,n+1} \\
z_{n+1, 1}& \cdots & z_{n+1, n} & z_{n+1, n+1}
\end{pmatrix}.
$$

The divisor $\{z_{n+1, n+1} =0\}$ in $\overline{\0}_\textup{min}$ has two components: $D_1$ (the last row is zero) and $D_2$ (the last column is zero). One can easily check that:
\begin{enumerate}
\item The map $\nu$ is an isomorphism outside $D_1 \cup D_2$.

\item The image $\nu(D_1)$ is the singular locus $S$ of $\cZ_X$. The fibers of $\nu|_{D_1}\colon  D_1 \to S$ are $\bA^1$ outside the origin and $\bA^n$ over the origin.

\item The map $\nu|_{D_2}\colon  D_2 \to \nu(D_2)$ is an isomorphism.

\item The map $\nu$ is not surjective. The closure of $\cZ_X \setminus {\rm Im}(\nu)$ is a divisor $W$ consisting of matrices with the last row  being identically zero. Moreover, we have $\cZ_X \setminus {\rm Im}(\nu) = W \setminus S$.
\end{enumerate}

As the matrix $A$ is unimodular, we have a symplectic resolution given by $\cY(A, \theta,0) \to \cY(A, 0, 0)$ for a general $\theta$.
Write $A=[\textbf{a}_1,\dots,\textbf{a}_{n+2}]$ into columns, namely
$$
\textbf{a}_1=\cdots=\textbf{a}_n= \begin{pmatrix} 1 \\ -1 \end{pmatrix}, \textbf{a}_{n+1} = \begin{pmatrix} 1 \\ 0 \end{pmatrix}, \textbf{a}_{n+2} = \begin{pmatrix} 0 \\ 1 \end{pmatrix}.
$$
The $\theta$-semi-stable set is given by
$$
(\mathbb{C}^{2n+4})^{\theta-ss} = \left\{ (\textbf{x}, \textbf{y}) \in \mathbb{C}^{2n+4}\,\middle\vert\, \theta \in \sum_{i: x_i \neq 0} \mathbb{Q}_{\geq 0} \textbf{a}_i + \sum_{j: y_j \neq 0} \mathbb{Q}_{\geq 0} (-\textbf{a}_j) \right\}. 
$$
Choose $\theta = \begin{pmatrix} 1 \\ -1 \end{pmatrix}$, then we have
$$
\Phi^{-1}(0)^\theta = \left\{ (\textbf{x}, \textbf{y}) \in \mathbb{C}^{2n+4} \,\middle\vert\, (x_1, \cdots, x_{n+1}) \neq 0, y_{n+2} \neq 0 \right\} \cap \Phi^{-1}(0).
$$
Then $\cY(A, \theta, 0) = \Phi^{-1}(0)^\theta /\!\!/\mathbb{T}^2$ and the symplectic resolution $\beta: \cY(A, \theta, 0) \to \cY(A, 0,0)$ is induced by
$$
\Phi^{-1}(0)^\theta \ni (\textbf{x}, \textbf{y}) \mapsto \begin{pmatrix} 
x_1y_1 & \cdots & x_1y_n & -x_1y_{n+1}x_{n+2} \\
\vdots & \vdots  & \vdots &  \vdots \\
x_ny_1 & \cdots & x_ny_n & -x_ny_{n+1}x_{n+2} \\
x_{n+1}y_1y_{n+2}& \cdots & x_{n+1}y_{n}y_{n+2} & x_{n+1}y_{n+1}
\end{pmatrix}.
$$
Now we define an injective morphism $\mathcal{H} \to \Phi^{-1}(0)^\theta$ as follows:
$$
\mathbb{C}^{2n+2} \supset \cH \ni (\textbf{x}, \textbf{y}) \mapsto (x_1, \cdots, x_{n+1}, -x_{n+1}y_{n+1}, y_1, \cdots, y_{n+1}, 1).
$$
This map descends to the quotients to give an injective morphism $$\tilde{\nu}:  T^*{\bP^n} \to \mathcal{Y}(A,\theta,0),$$
which makes the following diagram commutative:

\[
\begin{tikzcd}[column sep=large, row sep=large]
    T^*{\bP^n}  \arrow[r,"\widetilde{\nu}"] \arrow[d,"\Phi_{\bP^n}" left]
        &  \mathcal{Y}(A,\theta,0) \arrow[d,"\beta"] \\
    \overline{\0}_\textup{min} \arrow[r, "\nu"]
        & \cZ_X \simeq \cY(A,0,0)
\end{tikzcd}
\]

Now we compute the central fiber of $\beta$.  Assume $(\textbf{x}, \textbf{y}) \in \Phi^{-1}(0)^\theta$ which is mapped to 0. There are two cases:

\textbf{Case (i)}: \emph{There exists an index $1 \leq i \leq n$ such that $x_i \neq 0$.} Then we get $y_1=\cdots=y_n=0$ and $x_{n+1}y_{n+1} = x_{n+2}y_{n+2}=y_{n+1}x_{n+2}=0$.  As $y_{n+2} \neq 0$, we obtain two components in this case, both are isomorphic to $\mathbb{P}^n$:
$$
E_1=\{x_{n+2}=0, y_1=\cdots=y_n=y_{n+1}=0\}/\!\!/\mathbb{T}^2,
$$
which is in fact the image of the zero-section $\mathbf{0}_{\bP^n}$, and  
$$
E_2=\{x_{n+1}=x_{n+2}=0, y_1=\cdots=y_n=0\}/\!\!/\mathbb{T}^2,
$$
which is in fact the compactification of the $\mathbb{A}^n$-fiber of  $\nu|_{D_1}$ over the origin.

\textbf{Case (ii)}: \emph{$x_1=\cdots=x_n=0$.} Then we have $x_{n+1} \neq 0$ (as $(x_1, \cdots, x_{n+1}) \neq 0).$
Recall that $y_{n+2} \neq 0$, we must have $y_1=\cdots = y_n=y_{n+1}=0$ and $x_{n+2}=0$. But this subset is already contained in case (i). 

In conclusion,  the central fiber of $\beta$ consists of two $\bP^n$'s. {  For a non-zero singular point $s \in S={\rm Sing}(\cZ_X)$, $\cZ_X$ has a double point of type $A_1$ at $s$ by for example \cite[Theorem 3.6]{Nag21}, hence the fiber of $\beta$ over $s$ is just $\bP^1$.}
Now, we can obtain two more symplectic resolutions by performing Mukai flops along the two central $\bP^n$'s. These two symplectic resolutions correspond to  different choices of $\theta$. This gives three symplectic resolutions for $\cZ_X$.

In general, if $X$ is the blow-up of $n+1-k$ toric fixed points of $\bP^n$ (with $k \leq n+1$), then 
the symplectic singularity $\cZ_X$ is then isomorphic to the following
\[
\left\{ (z_{ij})_{1\leq i,j\leq n+1}\,\left\vert\, \sum_{i=1}^{n+1} z_{ii} =0, \quad 
{\rm rk} \begin{pmatrix} 
z_{1,1} & \cdots & \cdots  &\cdots &  \cdots\\
\vdots & \vdots  & \vdots &  \vdots & \vdots\\
z_{k,1} & \cdots & z_{k,k} & \cdots &  \cdots \\
z_{k+1,1} & \cdots   & \cdots & -z^2_{k+1,k+1} & \cdots   \\
\vdots & \vdots  & \vdots &  \vdots & \vdots\\
z_{n+1, 1}& \cdots & \cdots &  \cdots &-z_{n+1, n+1}^2
\end{pmatrix}  \leq 1 \right\}\right..
\]
We leave it to the readers to describe the symplectic resolutions in this general case.

\section{Symplectic singularities from horospherical varieties}
\label{s.horospherical-var}

We study $\cZ_X$ for smooth projective horospherical $G$-varieties $X$. We refer the reader to \cite[\S\,7]{Timashev2011} for general facts about horospherical varieties.

\subsection{Horospherical varieties and their moment maps}
\label{ss.def-horospherical}

We briefly recall the definition of horospherical varieties and some basic constructions used in our proof, see \cite[\S\,8.5]{Timashev2011}.

Let $G$ be a reductive and connected linear algebraic group and let $B$ be a Borel subgroup of $G$. A closed subgroup $H$ of $G$ is said to be \emph{horospherical} if it contains the unipotent radical of a Borel subgroup of $G$ and we also say that the homogeneous space
$G/H$ is horospherical. Denote by $P$ the normalizer $N_G(H)$ of $H$ in $G$. Then $P$ is a parabolic subgroup of $G$ such that $P/H$ is a torus. Thus $G/H$ is a torus bundle over the flag variety $W\coloneqq G/P$. The dimension $r$ of $P/H$ is called the \emph{rank} of $G/H$. A normal variety $X$ with an action of $G$ is said to be a \emph{horospherical $G$-variety} if $G$ has an open orbit isomorphic to $G/H$ for some horospherical subgroup $H$. { The torus $\bT^r\coloneqq P/H$ acts on the open orbit $G/H$ by right multiplication, so it acts freely and transitively on the fibers of $G/H\rightarrow G/P$.}

\begin{prop}[\protect{\cite[Proposition 4.1]{Pasquier2016}}]
\label{p.struc.horosspherical}
Let $X$ be a horospherical $G$-variety with open orbit $X_o=G/H$ isomorphic to a torus fibration $X_o\rightarrow G/P$ with $P$ a parabolic subgroup of $G$. Then there exists a smooth toric $\bT^r$-variety $Y$ and a
$G$-equivariant birational morphism $f$ from the smooth horospherical $G$-variety $G\times_P Y$ to $X$.
\end{prop}

Let $\fg$, $\fh$ and $\fp$ be the Lie algebras of $G$, $H$ and $P$, respectively. Then the cotangent bundle $T^*X_o$ can be canonically identified with $G\times_H \fh^{\perp}$. The moment map $T^*X_o\rightarrow \fg^*$ can be written as
\[
\Phi_{X_o}^G\colon T^*X_o=G\times_H \fh^{\perp} \longrightarrow \fg^*,\quad (g,\alpha)\longmapsto \Ad_g(\alpha)
\]
and it factors as
\begin{equation}
\label{e.moment-map-horospherical}
    \Phi_{X_o}^G\colon T^*X_o = G\times_H \fh^{\perp} \stackrel{q}{\longrightarrow} G\times_P \fh^{\perp} \stackrel{\psi}{\longrightarrow} \cM_X^G\subset \fg^*.
\end{equation}
The morphism $\psi$ is known to be generically finite and proper \cite{Richardson1974} (see also \cite[Theorem 8.7]{Timashev2011}). We denote by $G\times_P \fh^{\perp} \rightarrow \widetilde{\cM}_X\rightarrow \cM_X^G$ the Stein factorization of $\psi$. Denote the homogeneous vector bundle $G\times_P \fh^{\perp}$ over $W$ by $E$. Then we have $\widetilde{\cM}_X\cong \Spec R(\sO_{\bP E^*}(1))$.

\subsection{Cox rings}
\label{ss.Proof-Horo}
{ Now we aim to describe the Cox rings of $X$ and $\bP T_X$, using Theorem \ref{t.HS-Coxring} and the constructions from the previous subsection.} Let $X$ be an $n$-dimensional smooth projective $G$-horospherical variety with the open orbit $X_o=G/H$ and $X_o\rightarrow W\coloneqq G/P$ the induced torus bundle. 

\subsubsection{The Cox ring of $X$}

The following result can be directly derived from \cite[Theorem 4.3.2]{Brion2007} and \cite[Theorem 3.8]{Gagliardi2014}. We give an alternative proof by using Theorem \ref{t.HS-Coxring}. Recall that \emph{a boundary divisor} on a horospherical $G$-variety is a divisor whose support is contained in the boundary.
 
\begin{prop}
\label{p.Fanotype-horospherical}
    Let $X$ be a smooth projective $G$-horospherical variety. Then we have
    \[
    \Cox(X)\cong \Cox(W)[S_1,\dots,S_m],
    \]
    where  $S_{i}$ corresponds to the canonical section $1_{D_i}$ for each irreducible boundary divisor $D_i$ of $X$. In particular, the variety $X$ is of Fano type.
\end{prop}

\begin{proof}
    By Proposition \ref{p.struc.horosspherical}, the torus $\bT^r$ acts on $X$ with $X^{\circ}=X_o$ (see Theorem \ref{t.HS-Coxring}) and $X_o/\bT^r=W$. Thus the description of $\Cox(X)$ follows from Theorem \ref{t.HS-Coxring}. On the other hand, since $W$ is a Fano manifold, the affine variety $\Spec \Cox(W)$ has klt singularities by Theorem \ref{t.GOST-Cox-Fano} so that the affine variety
    \[
    \Spec \Cox(X) \cong \Spec \Cox (W)\times \bA^m
    \]
    also has only klt singularities, hence $X$ is of Fano type by Theorem \ref{t.GOST-Cox-Fano}.
\end{proof}

\begin{rem}
    It was proved by A.~Vezier in \cite[Corollary 2.4.2]{Vezier2023} that all the $\bQ$-factorial projective spherical varieties are of Fano type.
\end{rem}

\subsubsection{The Cox ring of $\bP T_X$}

We will relate the Cox ring of $\bP T_X$  to $\Cox(\bP E^*)$.

\begin{lem}
\label{l.FanoPE}
    The dual bundle $E^*$ is big and globally generated. In particular, the projective bundle $\bP E^*$ is a weak Fano manifold, i.e., $-K_{\bP E^*}$ is big and nef.
\end{lem}

\begin{proof}
    Note that $\fp^{\perp}\subset \fh^{\perp}$ is an invariant subspace under $P$ and $P$ acts trivially on the $r$-dimensional quotient $\fh^{\perp}/\fp^{\perp}$. Therefore we have an exact sequence
    \[
    0\longrightarrow T^*_{W} \cong G\times_P \fp^{\perp} \longrightarrow G\times_P \fh^{\perp}=E \longrightarrow G\times_P (\fh^{\perp}/\fp^{\perp})\cong \sO_{W}^{\oplus r} \longrightarrow 0.
    \]
    Taking dual yields
    \begin{equation}
    \label{e.E-exact-seq}
        0\longrightarrow \sO_{W}^{\oplus r} \longrightarrow E^* \longrightarrow T_{W} \longrightarrow 0.
    \end{equation}
    As $H^1(W,\sO_{W})=0$ and $T_{W}$ is globally generated, the vector bundle $E^*$ is also globally generated. 

    On the other hand, by \eqref{e.E-exact-seq}, the subvariety $\bP T_{W} \subset \bP E^*$ is a complete intersection of $r$ general elements in the sub-linear system of $|\sO_{\bP E^*}(1)|$ corresponding to the subspace
    \[
    H^0(W,\sO_{W}^{\oplus r})\subset H^0(W,E^*)\cong H^0(\bP E^*,\sO_{\bP E^*}(1)).
    \]
    So we have
    \[
    \deg(\sO_{\bP E^*}(1)) = \deg (\sO_{\bP E^*}(1)|_{\bP T_{W}}) = \deg(\sO_{\bP T_{W}}(1)).
    \]
    since $T_{W}$ is big and nef, we get $\deg(\sO_{\bP E^*}(1))>0$. Hence $\sO_{\bP E^*}(1)$ is big and $\bP E^*$ is a weak Fano manifold as $\sO_{\bP E^*}(-K_{\bP E^*})\cong \sO_{\bP E^*}(n)$ (see \cite[\S\,7.3.A]{Lazarsfeld2004a}).
\end{proof}

\begin{lem}
    Let $D$ be a prime divisor on $\bP T_X$. Then the generic isotropy group of $D$ under the $\bT^r$-action is trivial.
\end{lem}

\begin{proof}
    Note that the map $f$ in Proposition \ref{p.struc.horosspherical} is $\bT^r$-equivariant. Thus we may replace $X$ by $G\times_P Y$. Moreover, since the statement is clear if $D\cap \bP T_{X_o}\not=\emptyset$, we may assume that $D'\coloneqq \widebar{\pi}(D)$ is an irreducible boundary divisor contained in $X\setminus X_o$. Then there exists an irreducible toric divisor $F$ in $Y$ such that $D'=G\times_P F$ and the isotropy subgroup $\bT_{x}$ of the $\bT^r$-action at a general point $x\in D'$ is a one-dimensional torus as $F$ is of codimension one. Moreover, the induced action of $\bT_x$ on $T_x X$ yields a decomposition $T_{D',x}\oplus \bC$, where $\bT_{x}$ acts trivially on the first factor $T_{D',x}$ and acts with weight $\in\{-1,1\}$ on the second factor. So the isotropy subgroup of a point $y\in \bP T_{X,x}\setminus \bP T_{D',x}\subset D$ is trivial and we are done. 
\end{proof}

\begin{lem}
\label{l.separationPE}
    The variety $\bP E^*$ is a separation of the orbit space of the $\bT^r$-action on $\bP T_X$.
\end{lem}

\begin{proof}
    Note that the map $f$ in Proposition \ref{p.struc.horosspherical} is $\bT^r$-equivariant. Thus we may replace $X$ by $G\times_P Y$ and consider the rational map $g\colon \bP T_X\dashrightarrow \bP E^*$ induced by the projectivization $\bP T_{X_o}\rightarrow \bP E^*$ of $q\colon T^*X_o\rightarrow E$. It is enough to show that for { any prime $G$-stable divisor $D$} contained in $\bP T_X\setminus \bP T_{X_o}$, the general fiber of $D\dashrightarrow g(D)$ is exactly a $\bT^r$-orbit. Let $F$ be the irreducible toric divisor in $Y$ such that $D'=G\times_P F$, where $D'=\widebar{\pi}(D)$ is a prime divisor in $X$. 
    
    Let $\bT_{x} (\cong \bT)$ be the generic isotropy subgroup of the $\bT^r$-action on $D'$ and let $x\in D'$ be a general point contained in the fiber of $D'=G\times_P F\rightarrow W$ over the point $[P]\in W$. Then the isotropy subgroup $G_{x}$ of the $G$-action on $D'$ at $x$ is the kernel of the following composition
    \[
    P\longrightarrow P/H \cong \bT^r \longrightarrow \bT^r/\bT_{x}.
    \]
    So $P\subset N_{G}(G_x)$ and $H\subset G_x$. Denote by $D'_o\subset D'$ the unique $G$-open orbit. Then the $\bT_x$-action on $T^*X|_{D'_o}$ induces a decomposition $T_X|_{D'_o}=T_{D'_o}\oplus \bC_{D'_o}$, where $\bC_{D'_o}$ is the trivial bundle of rank one, such that $\bT_x$ acts trivially on the first factor and acts with weight $\in \{1,-1\}$ on the second factor. In particular, the general fiber of the natural projection
    \[
    (T^*X\setminus \textbf{0}_X)|_{D'_o} \longrightarrow T^*D'_o
    \]
    is exactly a $\bT_x$-orbit. Note that the restriction of the moment map $\varphi_X$ to $D'_o$ factors through the following map
    \[
    T^*D'_o\cong G\times_{G_x} \fg_x^{\perp} \longrightarrow G\times_P \fg_x^{\perp} \subset G\times_P \fh^{\perp}=E
    \]
    whose fibers are exactly the $(\bT^r/\bT_x)$-orbits. As a consequence, the fiber of the composition 
    \[
    \varphi_X\colon (T^*X\setminus \textbf{0}_X)|_{D'_o} \longrightarrow T^*D'_o \longrightarrow G\times_P \fg_x^{\perp}\subset E
    \]
    is a $\bT^r$-orbit and hence the general fiber of $D\dashrightarrow g(D)$ is also a $\bT^r$-orbit.
\end{proof}

Now we are in the position to finish the proof of Theorem \ref{t.Fanotype-MDS-horospherical}.

\begin{proof}[Proof of Theorem \ref{t.Fanotype-MDS-horospherical}]
    By Proposition \ref{p.Fanotype-horospherical}, it only remains to consider $\bP T_X$. By Theorem \ref{t.GOST-Cox-Fano} and Lemma \ref{l.FanoPE}, the Cox ring $\Cox (\bP E^*)$ is finitely generated and $\Spec\Cox (\bP E^*)$ has only klt singularities. In particular, by Theorem \ref{t.HS-Coxring} and Lemma \ref{l.separationPE}, the Cox ring $\Cox (\bP T_X)$ is also finitely generated and hence $\bP T_X$ is a Mori dream space. If additionally $X$ does not contain any $G$-stable divisor, namely $\codim(X\setminus X_o)\geq 2$, then Theorem \ref{t.HS-Coxring} implies that $\Cox (\bP T_X)\cong \Cox (\bP E^*)$ and hence $\bP T_X$ is of Fano type by Theorem \ref{t.GOST-Cox-Fano}. 
\end{proof}

In general it is difficult to describe the geometry of $\cZ_X$ explicitly. However, since $\cZ_X$ admits a natural $\bT^r$-action, we have the following description of the quotient of $\cZ_X$ by $\bT^r$.

\begin{prop}
\label{p.QuotientZX}
    Let $X$ be a smooth projective $G$-horospherical variety without $G$-stable divisors. Then $\widetilde{\cM}_X\cong \cZ_X/\!\!/\bT^r$.
\end{prop}

\begin{proof}
    It remains to prove $S(X)^{\bT^r}\cong R(\sO_{\bP E^*}(1))$. As $X$ does not contain $G$-stable divisors, we have $\codim(X\setminus X_o)\geq 2$ and thus the natural restriction $S(X)\rightarrow S(X_o)$ is an isomorphism of graded $\bC$-algebras. As $T^*X_o/\bT^r=E$, we obtain
    \[
    S(X_o)^{\bT^r}\cong H^0(T^*X_o,\sO_{T^*X_o})^{\bT^r}\cong H^0(E,\sO_E)\cong R(\sO_{\bP E^*}(1)),
    \]
    where the second isomorphism is given by the pull back of $q\colon T^*X_o\rightarrow E$.
\end{proof}

\begin{question}
    The following two questions are inspired by Question \ref{q.Main-question} and Theorem \ref{t.Fanotype-MDS-horospherical}.
    \begin{enumerate}
        \item Is $\bP T_X$ of Fano type for \emph{any} smooth projective horospherical variety $X$?

        \item Find out all smooth projective horospherical varieties $X$ such that $\cZ_X$ admits a symplectic resolution.
    \end{enumerate}
\end{question}

\subsection{Picard number one case}
\label{ss.Proof-Picardone}

In the following we focus on smooth projective $G$-horospherical varietes with Picard number one, which are completely classified by B.~Pasquier. Without loss of generality, we shall assume that $X$ is not homogeneous and then there are exactly five different types (\cite[Theorem 0.1]{Pasquier2009}):
\begin{itemize}
    \item $(B_m,\omega_{m-1},\omega_m)$ $(m\geq 3)$;

    \item $(B_3,\omega_1,\omega_3)$;

    \item $(C_m,\omega_i,\omega_{i+1})$ $(m\geq 2, 1\leq i\leq m-1)$;

    \item $(F_4,\omega_2,\omega_3)$;

    \item $(G_2,\omega_2,\omega_1)$.
\end{itemize}

We refer the reader to \cite{Pasquier2009} for the notations and we give a brief introduction to the geometry of $X$. There are exactly three $G$-orbits in $X$: the unique open $G$-orbit $X_o=G/H$ and two closed $G$-orbits $Z=G/P_Z$ and $Y=G/P_Y$, where $P_Y, P_Z$ are two parabolic subgroups of $G$ such that
{  $Z$ is ${\rm Aut}(X)$-stable while $Y$ is not (recall that $X$ has only two ${\rm Aut}(X)$-orbits),} Moreover,  both $Z$ and $Y$ have codimension at least two in $X$, so $\cZ_X$ is symplectic by Theorem \ref{t.criterion-contact-symplectic} and Theorem \ref{t.Fanotype-MDS-horospherical}. The open subset $X_Y\coloneqq X\setminus Z$ is isomorphic to a homogeneous vector bundle 
\[
\phi_Y\colon X_Y\cong G\times_{P_Y} V_Y\longrightarrow Y=G/P_Y,
\]
where $V_Y$ is a simple $P_Y$-module, $Y\subset X_Y$ is just the zero section of $\phi_Y$ and the open subset $X_o$ can be identified to $X_Y\setminus Y$, which satisfies the following commutative diagram:
\[
\begin{tikzcd}[column sep=large, row sep=large]
    X_o=G/H \arrow[r,hookrightarrow] \arrow[d]
        &  X_Y\cong G\times_{P_Y} V_Y \arrow[d,"\phi_Y"] \\
    W\coloneqq G/P \arrow{r}
        & Y=G/P_Y
\end{tikzcd}
\]
The torus $\bT\coloneqq P/H$ has rank one and it acts on the fibers of $\phi_Y$ by scaling. Since $\phi_Y\colon X_Y\rightarrow Y$ can be naturally identified to the normal bundle of $Y$, the vector group $\bV\coloneqq H^0(Y,N_{Y/X})$ is exactly the space of sections of $\phi_Y$ and  it acts naturally on $X_Y$ as follows:
\[
s\cdot (g,v) \longmapsto (g,v+s(g)),
\]
where $s\in \bV$. This action extends to the whole $X$ as $\codim Z\geq 2$. Denote by $\widetilde{G}$ the automorphism group of $X$ with $\widetilde{\fg}$ its Lie algebra. By \cite[Theorem 1.11]{Pasquier2009}, we have a natural isomorphism 
\[
\widetilde{\fg} \cong \fg \oplus \ft \oplus \fv,
\]
where $\ft$ and $\fv$ are the Lie algebras of $\bT$ and $\bV$, respectively. Let $\widebar{G}\coloneqq G\times \bT$. Then there exists a homomorphism $\widebar{G} \rightarrow \widetilde{G}$ with finite kernel. Then $X_o$ is also the unique open $\widebar{G}$-orbit. Let $\widebar{H}$ be the isotropic subgroup of $\widebar{G}$ at $[H]\in X_o$. Then $P\subset N_{\widebar{G}}(\widebar{P})$ and the natural projection to the $\fg$-factor induces an isomorphism $\widebar{\fh}^{\perp} \rightarrow \fh^{\perp}$ such that 
\[
T^*X_o \cong \widebar{G}\times_{\widebar{H}} \widebar{\fh}^{\perp} \cong G\times_H \widebar{\fh}^{\perp}.
\]
So the moment map $\Phi_X^{\widebar{G}}\colon T^*X_o \rightarrow \cM_X^{\widebar{G}}$ factors as
\[
T^*X_o = G\times_H \widebar{\fh}^{\perp} \longrightarrow E\cong G\times_P \widebar{\fh}^{\perp} \stackrel{h}{\longrightarrow} \cM_X^{\widebar{G}}\subset \widebar{\fg}^*.
\]
The map $h$ is generically finite and proper such that its Stein factorization coincides with $E\rightarrow \widetilde{\cM}_X$. So we have the following commutative diagram:
\[
\begin{tikzcd}[column sep=large, row sep=large]
    T^*X_o \arrow[r, phantom, "\subset"] \arrow[dd,"q" left]
        & T^*X \arrow[r,"\Phi_X^{\widetilde{G}}"] \arrow[dr,"\Phi_X^{\widebar{G}}"] \arrow[ddr,"\Phi_X^G" near start, swap]
            & \cM_X^{\widetilde{G}}\arrow[r, phantom,"\subset"] \arrow[d,"p_{\fv}"]
                & \widetilde{\fg}^*=\fg^*\oplus \ft^*\oplus \fv^* \arrow[d,"p_{\fv}"] \\
        &   & \cM_X^{\widebar{G}} \arrow[r, phantom, "\subset"]  \arrow[d,"p_{\ft}"] 
                & \widebar{\fg}^*=\fg^*\oplus \ft^* \arrow[d,"p_{\ft}"] \\
    E=G\times_P \fh^{\perp} \arrow[r,"\widetilde{\psi}" below] \arrow[rr, bend right, "\psi"] \arrow[urr,"h"]
        &\widetilde{\cM}_X \arrow[r,"g" below] \arrow[ur,"\widebar{g}" below]
            & \cM_X^G \arrow[r, phantom, "\subset"]
                & \fg^*.
\end{tikzcd}
\]

\begin{lem}
\label{l.genericfiniteness}
    The map $\Phi_X^{\widetilde{G}}$ is generically finite. Moreover, if $\Phi_{W}^G$ is birational onto its image, then the map $\Phi_{X}^{\widetilde{G}}$ is also birational.
\end{lem}

\begin{proof}
    As $h^{-1}(\cM_X^{\widebar{G}}\cap \fg^*)=T^*W$ and the restriction $h|_{T^*W}\colon T^*W\rightarrow \fg^*$ coincides with $\Phi_W^G$, it suffices to show that a general fiber of $q$ is mapped birationally onto its image in $\cM_X^{\widetilde{G}}$ under $\Phi_X^{\widetilde{G}}$. Fix a general fiber of $\phi_Y$, which is isomorphic to the vector space $V_Y$. Note that $\bT$ acts on $V_Y$ by multiplication and the evaluation map induces a quotient $\bV\rightarrow V_Y$ such that $\bV$ acts on $V_Y$ via the additive action of $V_Y$. In particular, we have the moment map of $\bV$ restricted to $V_Y$ as follows:
    \[
    \Phi_X^{\bV}|_{V_Y}\colon T^*X|_{V_Y} \stackrel{\nu}{\longrightarrow} T^*V_Y=V_Y\times V^*_Y \stackrel{p_2}{\longrightarrow} V_Y^*\subset \fv^*,
    \]
    where the projection $p_2$ coincides with the moment $\Phi_{V_Y}^{\bV}$ (cf. Example \ref{e.momentmap}). Let $O\subset T^*X|_{V_Y}$ be a general $\bT$-orbit. Since both $\nu$ and $p_2$ are $\bT$-equivariant, the image $p_2\circ \nu(O)$ is again a $\bT$-orbit and the induced morphism $O\rightarrow p_2\circ \nu(O)$ is an isomorphism.
\end{proof}

\begin{rem}
{  The degree of the Springer map $\Phi_W^G\colon T^*W\rightarrow \cM_W^G=\widebar{\cO}_P$ is given by Hesselink  (\cite[Theorem 7.1(d)]{Hesselink1978}) when $G$ is of classical type. For  $(F_4,\omega_2,\omega_3)$, the image $\widebar{\cO}_P$ is the nilpotent orbit closure labeled $F_4(a_2)$, which is an even orbit, hence the Springer map is birational. For $(G_2,\omega_2,\omega_1)$, $W$ is just $G/B$, hence the Springer map is birational.  These results are summarized in the following table.}

\footnotesize
\renewcommand*{\arraystretch}{1.6}
			\begin{longtable*}{|M{6cm}|M{3.5cm}|M{1cm}|M{1.5cm}|}
				\hline
				
				$X$
				& $m$
                & $i$
				& $\deg(\Phi_W^G)$
				\\
                \hline

                \multirow{2}{*}{$(B_m,\omega_{m-1},\omega_m) (m\geq 3)$}
				& \textup{even}
                & -
				& 1
				\\
				\cline{2-4}

				& \textup{odd}
                & -
				& 2
				\\
				\hline

                $(B_3,\omega_1,\omega_3)$
                & -
                & -
                & 2
                \\
                \hline

                \multirow{4}{*}{$(C_m,\omega_i,\omega_{i+1}) (m\geq 2, 1\leq i\leq m-1)$}
                & $i=m-1$
                & -
                & 1
                \\
                \cline{2-4}

                & $\frac{2m-2}{3}\leq i \leq m-2$
                & -
                & 2
                \\
                \cline{2-4}

                & \multirow{2}{*}{$i<\frac{2m-2}{3}$}
                & \textup{even}
                & 2
                \\
                \cline{3-4}

                &
                & \textup{odd}
                & 4
                \\
                \hline

                $(F_4,\omega_2,\omega_3)$
                & -
                & -
                & 1
                \\
                \hline

                $(G_2,\omega_2,\omega_1)$
                & -
                & -
                & 1
                \\
                \hline
			\end{longtable*}
\normalsize
\end{rem}

\begin{proof}[Proof of Proposition \ref{p.PicardnubmeroneHorospherical}]
    We have $\cM_X^{\widebar{G}}\cong \cM_X^{\widetilde{G}}/\!\!/\bT$ as $\widebar{\fg}=\widetilde{\fg}^{\bT}$. Then the result follows from Proposition \ref{p.QuotientZX} and Lemma \ref{l.genericfiniteness}.
\end{proof}

\section{Symplectic singularity from the quintic del Pezzo threefold}
\label{s.V5}

In this section we consider the quintic del Pezzo threefold $X\subset \bP^6$, i.e., a smooth codimension three linear section of $\Gr(2,5)\subset \bP^{9}$. {In \S\,\ref{ss.lines}, we recall some basic facts about the lines contained in $X$. In \S\,\ref{ss.weak-Fano-model}, we review the construction of weak Fano model of $\bP T_X$ from \cite{HoeringPeternell2024}. In \S\,\ref{ss.anti-canonical-model}, we study the geometry of the anti-canonical model of $\bP T_X$ via its weak Fano model. Finally, we combine the various results to describe the $\bQ$-factorial terminalization of $\cZ_X$ in \S\,\ref{ss.Q-fact-termi}.}

\subsection{Geometry of lines on $X$}
\label{ss.lines}

We briefly recall some facts about the lines on $X$ and refer the reader to \cite{FurushimaNakayama1989} for more details. We denote by $\sO_X(1)$ the restriction $\sO_{\bP^{6}}(1)|_{X}$, which is the ample generator of $\Pic(X)\cong \bZ$. The Fano variety of lines on $X$ is isomorphic to $\bP^2$ and the universal family $u\colon \cU\rightarrow \bP^2$ is a $\bP^1$-bundle such that the evaluation morphism $\nu\colon \cU\rightarrow X$ is a degree three finite morphism. Given a line $l$ in $X$,  the restriction $T_X|_l$ is isomorphic to one of the following:
\[
\sO_{\bP^1}(2)\oplus \sO_{\bP^1}\oplus \sO_{\bP^1} \quad\textup{or} \quad \sO_{\bP^1}(2)\oplus \sO_{\bP^1}(1)\oplus \sO_{\bP^1}(-1).
\]
A line $l$ is said to be \emph{of second type} if $T_X|_l$ is of the second type. The subset $B$ of $\bP^2$ parametrizing lines of second type is a smooth conic curve. 

The surface $R\coloneqq u^{-1}(B)$ is isomorphic to $B\times \bP^1$ so that $D=\nu(R)\in |\sO_X(2)|$ is exactly the branch locus of $\nu$ and $\nu^{-1}(D)$ decomposes into $R'\cup R$, where $R'\cong \bP^1\times \bP^1$, such that the induced morphism $R'\rightarrow \bP^2$ is a double covering. Consider the following exact sequence of vector bundles
\begin{equation}
\label{e.defn-totaldualVMRT}
    0\longrightarrow T_{\cU/\bP^2} \longrightarrow \nu^*T_X \longrightarrow Q\longrightarrow 0.
\end{equation}
Then $q\colon \bP Q\rightarrow \cU$ is a $\bP^1$-bundle and $\cB\coloneqq f(\bP Q)$ is the \emph{total dual VMRT} of $X$, where $f\colon \bP \nu^*T_X\rightarrow \bP T_X$ \cite[\S\,2.2]{HoeringLiuShao2022}. In summary, we have the following commutative diagram:
\begin{equation}
\begin{tikzcd}[row sep=large, column sep=large]
        & \bP Q \arrow[r,"f"] \arrow[d, phantom, sloped, "\subset"]
            & \cB  \arrow[d, phantom, sloped, "\subset"]
                & \\
        & \bP \nu^* T_X \arrow[d,"q"] \arrow[r,"f"]
            & \bP T_X \arrow[d,"\widebar{\pi}"]
                & \\
    R \arrow[r,"\subset"{sloped}] \arrow[d,"u" left] \arrow[rrr, controls={+(0.3,-2.2) and +(-1,-0.5)}, "\nu" near start]
        & \cU \arrow[r,"\nu"] \arrow[d,"u"]
            & X 
                & D \arrow[l, phantom, "\supset"] \\
    B \arrow[r,phantom,"\subset"]
        & \bP^2
            & R' \arrow[l,"u"] \arrow[ul,"\supset" {sloped}] \arrow[ur,"\nu" below]
                &
\end{tikzcd}
\end{equation}
Given  a line $[l]\in \bP^2\setminus B$, recall that a \emph{minimal section} $\widebar{l}$ over $l$ is a section of $\bP T_X|_l \rightarrow l$ corresponding to a trivial quotient $T_X|_l\rightarrow \sO_{\bP^1}$. {It is clear from the the construction of $\cB$ that minimal sections form a covering family of curves on $\cB$.}

\begin{lem}
\label{l.ampleness-restriction-TX}
    Let $C\subset X$ be an irreducible curve such that $C\not\subset D$. Then either $C$ is a line of the first type or $Q|_{\widetilde{C}}$ is ample for any component $\widetilde{C}$ of $\nu^{-1}(C)$. 
\end{lem}

\begin{proof}
    As $\widetilde{C}\not\subset R$, the map $T_{\cU}|_{\widetilde{C}}\rightarrow \nu^*T_X|_{\widetilde{C}}$ is generically surjective and so the induced map $u^* T_{\bP^2}|_{\widetilde{C}}\rightarrow Q|_{\widetilde{C}}$ is generically surjective. As $T_{\bP^2}$ is ample, either $\widetilde{C}$ is $u$-vertical or $Q|_{\widetilde{C}}$ is ample. { In the former case we have an exact sequence
    \[
    0\longrightarrow T_{\widetilde{C}}\cong \sO_{\sO_{\bP^1}}(2) \longrightarrow T_{\cU}|_C \longrightarrow u^*T_{\bP^2}|_{\widetilde{C}} \cong \sO_{\bP^1}^{\oplus 2} \longrightarrow 0.
    \]
    Since $T_{\cU}|_{\widebar{C}}\rightarrow \nu^*T_X|_{\widebar{C}}$ is generically surjective, it follows that $\nu^*T_X|_{\widetilde{C}}$ does not contain negative factors and hence $C$ is a line of the first type.}
\end{proof}

\subsection{Weak Fano model of $\bP T_X$} 
\label{ss.weak-Fano-model}

Our starting point is the following result obtained by A.~H\"oring and T.~Peternell in \cite{HoeringPeternell2024} recently.

\begin{thm}[\protect{\cite[Theorem 3.1]{HoeringPeternell2024}}]
    Let $X$ be the quintic del Pezzo threefold. Then there exists an anti-flip $\psi\colon \bP T_X \dashrightarrow \cY$ such that $\cY$ is a $\bQ$-factorial weak Fano variety with terminal singularities.
\end{thm}

This theorem was proved in \cite{HoeringPeternell2024} by explicitly constructing  the variety $\cY$. Let us briefly explain their construction (see also the proof of \cite[Proposition 4.1]{HoeringPeternell2024}). Given a line $l$ of second type in $X$, we denote by $\widetilde{l}\subset \bP T_X$ the unique negative section corresponding to the quotient $T_X|_l\rightarrow \sO_{\bP^1}(-1)$. By \cite[Lemma 3.4]{HoeringPeternell2024}, the locus
\[
S\coloneqq \bigcup_{[l]\in B} \widetilde{l}\quad\subset\quad \bP T_X
\]
is isomorphic to $\bP^1\times \bP^1$. { There exists a birational contraction $\varphi^+\colon \bP T_X\rightarrow \cM$ defined by the big and semi-ample line bundle $\sO_{\bP T_X}(1)\otimes \widebar{\pi}^*\sO_X(1)$, which contracts $S$ to $\bP^1$ so that the fibers are exactly those $\widetilde{l} $(\cite[(13) and Lemma 3.3]{HoeringPeternell2024})}.

Denote by $\mu_1\colon \cX_1\rightarrow \bP T_X$ the blow-up along $S$ with exceptional divisor $E_1$. Then there exists a section $S_1$ of the projective bundle $E_1\rightarrow S$ (\cite[(22) and Lemma 3.10]{HoeringPeternell2024}). Denote by $\mu_2\colon \cX_2\rightarrow \cX_1$ the blow-up along $S_1$ with exceptional divisor $E_2$. Thanks to \cite[Lemma 3.11, Proposition 3.12 and Proposition 3.13]{HoeringPeternell2024}, there exists a contraction $\nu_2\colon \cX_2\rightarrow \cY_1$ to a projective manifold $\cY_1$, which is the blow-up along a smooth center $W_1$ of codimension two with $E_2$ being its exceptional divisor. Let $\widetilde{E}_1$ be the strict transform of $E_1$ in $\cY_1$. Then there exists a $\bP^2$-bundle structure $\widetilde{E}_1\rightarrow S'\cong \bP^1\times \bP^1$ and \cite[Lemma 3.15 and Proposition 4.1]{HoeringPeternell2024} imply that it can be extended to a birational contraction $\nu_1\colon \cY_1\rightarrow \cY$ such that
\begin{itemize}
    \item $\cY$ has only $\bQ$-factorial terminal singularities, and 
    
    \item $-K_{\cY}$ is nef and big, and

    \item the exceptional locus of $\nu_1$ is $\widetilde{E}_1$.
\end{itemize}
Denote by $W'$ the image $\nu_1(W_1)$. By \cite[Lemma 3.3 and Proposition 3.16]{HoeringPeternell2024}, we have the following commutative diagram:
\begin{equation}
\label{eq.Weak-Fano-model}
    \begin{tikzcd}[column sep=large, row sep=large]
    S\arrow[r,phantom,"\subset"] \arrow[ddrr]
        & \bP T_X \arrow[rd,"\varphi^+" swap] \arrow[rr,"\psi", dashed]
            & 
                & \cY \arrow[ld,"\varphi^-"]
                    & S'\cup W' \arrow[l,phantom,"\supset"] \arrow[ddll]\\
        &
           & \cM
               &
                   & \\
        &
           & T\cong \bP^1 \arrow[u,phantom,sloped,"\subset"]
               &
                   &            
\end{tikzcd}
\end{equation}
such that 
\begin{itemize}
    \item the morphism $\varphi^+$ is a small birational contraction with exceptional locus $S$ such that $\varphi^{+}|_S\colon S\cong \bP^1\times \bP^1\rightarrow T=\bP^1$ is one of the natural projections such that its fibers are the negative sections $\widetilde{l}$ over the lines of second type and $K_{\bP T_X}$ is $\varphi^+$-ample, and 

    \item the morphism $\varphi^-$ is a small birational contraction with exceptional locus $S'\cup W'$ such that $\varphi^{-1}|_{S'}\colon S'=\bP^1\times \bP^1\rightarrow T=\bP^1$ is one of the natural projections, $\varphi^{-}|_{W'}\colon W'\rightarrow T$ is $\bP^2$-bundle and $-K_{\cY}$ is $\varphi^-$-ample. 
\end{itemize}

\subsection{Anti-canonical model of $\bP T_X$}
\label{ss.anti-canonical-model}

Now we aim to describe the anti-canonical model of $\bP T_X$. { Since $\psi$ is an isomorphism in codimension one, we can choose a Weil divisor $L$ on $\cY$ such that $\sO_{\bP T_X}(\psi^* L)\cong \sO_{\bP T_X}(1)$.} Then $-K_{\cY}=3L$. By the Basepoint Free Theorem, the $\bQ$-Cartier divisor $L$ is semi-ample, so it defines a birational contraction
\[
\Phi_{mL}\colon \cY \longrightarrow \cY^{\textup{can}}.
\]
onto a normal variety for some $m\gg 1$. Then $-K_{\cY^{\textup{can}}}$ is ample and $\Phi_{mL}$ is crepant. Moreover, since $\psi$ is an isomorphism in codimension one, one can easily derive
\[
\Proj R(-K_{\bP T_X})\cong \Proj R(-K_{\cY}) \cong \Proj R( L) = \Proj(\cC_L) \cong \cY^{\textup{can}}.
\]
Denote by $\cB'$ the strict transform of $\cB$ in $\cY$ and $H'$ the strict transform of $\widebar{\pi}^*A$ for some $A\in |\sO_X(1)|$. 

\begin{lem}
\label{l.stric-transf-totaldualVMRT}
    Let $\widebar{l}$ be a general minimal section over a general line $[l]\in \bP^2\setminus B$. Then $\widebar{l}\cap S=\emptyset$. In particular, the divisor $\cB'$ is contained in $\Exc(\Phi_{mL})$.
\end{lem}

\begin{proof}
    As $D\in |\sO_X(2)|$, the intersection $l\cap D$ consists of two points, namely $x_1$ and $x_2$. On the other hand, as $S\rightarrow D$ is one-to-one (\cite[Lemma 3.2]{HoeringPeternell2024}), there exists a unique point $\widebar{x}_i\in S$ such that $\widebar{\pi}(\widebar{x}_i)=x_i$ for each $i=1,2$. Note that the minimal sections over $l$ are exactly the fibers of $\bP(N_{l/T_X})\cong \bP^1\times l \rightarrow l$, so a general minimal section over $l$ is disjoint from the $\widebar{x}_i$'s and hence also $S$. Let $\widebar{l}'$ be the strict transform of $\widebar{l}$ in $\cY$. As $\widebar{l}\cap S=\emptyset$, we must have
    \[
    0=-K_{\bP T_X}\cdot \widebar{l} = -K_{\cY}\cdot \widebar{l}' = 3L\cdot \widebar{l}'.
    \]
    So $\widebar{l}'$ is contracted by $\Phi_{mL}$ and hence it is contained in $\Exc(\Phi_{mL})$. { In particular, since the minimal sections cover $\cB$}, the strict transforms of minimal sections cover $\cB'$, therefore $\cB'$ is contained in $\Exc(\Phi_{mL})$.
\end{proof}

\begin{lem}
\label{l.exceptional-locus}
    $\Exc(\Phi_{mL})=\cB'$.
\end{lem}

\begin{proof}
    By \cite[Theorem 5.4]{HoeringLiuShao2022}, we have $\cB\in |\sO_{\bP T_X}(3)\otimes \widebar{\pi}^*\sO_X(-1)|$. On the other hand, as $\rho(\cY)=2$, any irreducible curve $C$ contracted by $\Phi_{mL}$ is numerically proportional to $\widebar{l}'$ by Lemma \ref{l.stric-transf-totaldualVMRT}, where $\widebar{l}'$ is the strict transform of a general minimal section $\widebar{l}$ over a general line $[l]\in \bP^2\setminus B$. Then we have
    \[
    \cB'\cdot C = a\cB'\cdot \widebar{l}' = a\cB\cdot \widebar{l} = - a A\cdot l < 0,
    \]
    where $a\in \bQ_{>0}$. So $C$ is contained $\cB'$ and we are done.
\end{proof}

\begin{lem}
\label{l.dimension-image-cC'}
    $\dim(\Phi_{mL}(\cB'))=3$.
\end{lem}

\begin{proof}
    Let $F$ be an irreducible component of a general fiber of $\cB'\rightarrow \Phi_{mL}(\cB')$. We want to show that $F$ is actually a curve. Assume to the contrary that $\dim (F)\geq 2$.

    According to the Step 2 and Step 3 in the proof of \cite[Proposition 4.1]{HoeringPeternell2024}, the restrictions $-K_{\cY}|_{S'}$ and $-K_{\cY}|_{W'}$ are ample. So $\dim (F\cap \Exc(\varphi^{-}))\leq 0$. { Let $C\subset F$ be a curve arising as a general complete intersection of very ample divisors on $F$ so that} $C$ is disjoint from $\Exc(\varphi^{-})$. By abuse of notation, we denote its strict transform in $\bP T_X$ again by $C$. Then $C\subset \cB$ is disjoint from $S=\Exc(\varphi^+)$ and satisfies 
    \[
    c_1(\sO_{\bP T_X}(1))\cdot C=L\cdot C=0.
    \]
    So $C$ is $\widebar{\pi}$-horizontal. Denote $\widebar{\pi}(C)$ by $C'$. Since $\cB\rightarrow X$ is surjective and $F$ is general, we can assume that { $C'\not \subset D$, where $D$ is the branch locus of the evaluation morphism $\nu\colon \cU\rightarrow X$}. 
    
    Let $\widetilde{C}'$ be an irreducible component of $\nu^{-1}(C')\subset \cU$. Then there exists an irreducible curve $\widetilde{C}\subset \bP Q$ such that $f(\widetilde{C})=C$ and $q(\widetilde{C})=\widetilde{C}'$. As $f^*\sO_{\bP T_X}(1)|_{\bP Q}\cong \sO_{\bP Q}(1)$, we obtain 
   \[
   c_1(\sO_{\bP Q}(1))\cdot \widetilde{C} = 0,
   \]
    so $Q|_{\widetilde{C}'}$ is not ample. Then Lemma \ref{l.ampleness-restriction-TX} implies that $C'$ is a line of the first type and consequently $C$ is a minimal section over $C'$.

    { Finally, for a general point $\widebar{x}\in \cB$ with $x=\widebar{\pi}(\widebar{x})\in X$, exactly three lines pass through $x$ since $\nu\colon \cU\rightarrow X$ has degree three. For any such line $l$ through $x$, there is at most one minimal section $\widebar{l}$ over $l$ containing $\widebar{x}$, implying at most three minimal sections through $\widebar{x}$ in total. This contradicts the assumption $\dim(F)\geq 2$, because the complete intersection curves through $\widebar{x}$ form a positive-dimensional family in this situation.}
\end{proof}

We recall that $\textbf{P} E\coloneqq \bP E^*$ for a vector bundle or a vector space $E$. Let $l$ be a line of the first type. Then $N_{l/X}\cong \sO_{\bP^1}\oplus \sO_{\bP^1}$ and we have the following composition of projection maps
\[
\Psi_l\colon \bfP T_X|_l \dashrightarrow \bfP N_{l/X}\cong l\times \bfP V \rightarrow \bfP V\cong \bP^1,
\]
where the first map is the projection and $V$ is a two-dimensional vector space. { Moreover, the tangent bundle $T_l\cong \sO_{\bP^1}(2)$ of $l$ is a direct summand of $T_X|_l$, which yields a section $\bfP T_l$ of $\bfP T_X|_l\rightarrow l$.}

\begin{lem}[\protect{\cite[\S\,5.3.1]{HoeringLiu2023}}]
\label{l.tangentdirections-lines-V5}
    Let $l$ be a general line of the first type, and let $\cC_{l}\subset \bfP T_X|_l$ be the tangent direction along $l$ of the lines meeting $l$. Then $\cC_{l}$ is a rational curve, which is disjoint from $\bfP T_l$, such that $\Psi_l\colon \cC_l\rightarrow \bfP V$ is an isomorphism. 

     { In other words, given an arbitrary one-dimensional subspace $V'$ of $V$, there exists a unique point $x\in l$ and a unique line $l'$ meeting $l$ at $x$ such that $V'$ is exactly  the image of $T_{l',x}$ under the natural quotient
    \[
    T_{X,x} \longrightarrow N_{l/X,x} = V.
    \]
    }
\end{lem}

{ Recall that the minimal sections over lines of the first type cover $\cB$  and their strict transforms cover $\cB'$. By the proof of Lemma \ref{l.stric-transf-totaldualVMRT}, the strict transform of a general minimal section of a general line is contracted by $\Phi_{mL}$. In the following we prove the converse.}

\begin{prop}
\label{p.fibres-C'}
    The general fiber of $\cB'\rightarrow \Phi_{mL}(\cB')$ is a union of two smooth rational curves meeting at a point, which are the strict transforms of minimal sections.
\end{prop}

\begin{proof}
    By Lemma \ref{l.dimension-image-cC'}, the general fiber $F$ of $h\colon \cB'\rightarrow \Phi_{mL}(\cB')$ is one-dimensional, so the irreducible components of $F$ are exactly strict transforms of minimal sections. In particular, since $F$ is connected, it suffices to show that for for a general line $[l]\in \bP^2\setminus B$ and a general minimal section $\widebar{l}$ over $l$, there exists a unique line $[l']\in \bP^2$ and a unique minimal section $\widebar{l}'$ over $l'$ such that $\widebar{l}'\cap \widebar{l}\not=\emptyset$.
    
    Recall that $\widebar{l}$ is a fiber of 
    \[
    \bP N_{l/X}\cong \bfP N^*_{l/X} = l\times \bfP V^*\longrightarrow \bfP V^*, 
    \]
    which corresponds a one-dimensional subspace $K_{\widebar{l}}$ of $V^*$, i.e., $\widebar{l}=l\times \{[K_{\widebar{l}}]\}$.
    
    { Let $l'$ be another line meeting $l$ at $x\in l$ such that there exists a minimal section $\widebar{l}'$ over $l'$ with 
    \[
    \emptyset\not=\widebar{l}\cap \widebar{l'}\in \bP T_{X,x}=\bfP T_{X,x}^*.
    \]
    By the definition of minimal sections, the point $\widebar{l'}\cap \widebar{l}$ corresponds to the subspace 
    \[
    K_x\coloneqq N^*_{l/X,x}\cap N_{l'/X,x}^* \subset T_{X,x}^*.
    \]
    In other words, under the canonical identification $N^*_{l/X,x}=V^*$, we have $K_x=K_{\widebar{l}}$ so that 
    \[
    \widebar{l}\cap \widebar{l}' = [K_x]=[K_{\widebar{l}}]\in \bfP(N^*_{l/X,x}) \subset \bfP T^*_{X,x} = \bP T_{X,x}\subset \bP T_X.
    \]

    Consequently, the image of $T_{l',x}$ in $N_{l/X,x}$ under the quotient $T_{X,x}\rightarrow N_{l/X,x}$ is exactly the kernel of the quotient map 
    \[
    N_{l/X,x}=V\longrightarrow K_x^* = K_{\widebar{l}},
    \]
    which is uniquely determined by $\widebar{l}$. Then Lemma \ref{l.tangentdirections-lines-V5} says that such $\widebar{l}'$ exists and is unique.}
\end{proof}

\subsection{$\bQ$-factorial terminalization of $\cZ_X$}
\label{ss.Q-fact-termi}
In this subsection we want to describe the singular locus of $\cZ_X$ and also its $\bQ$-factorial terminalizations. Recall that we have the following commutative diagram:
\[
\begin{tikzcd}[row sep=large, column sep=large]
    p^{-1}(\cB') \arrow[r,phantom,"\subset"] \arrow[d]
        & \cL \arrow[r,"\varphi_L"] \arrow[d,"p"]
            & \cC_L=\cZ_X \arrow[d, dashed] 
                & Z=\varphi_L(p^{-1}(\cB')) \arrow[d,dashed] \arrow[l,phantom,"\supset"] \\
    \cB' \arrow[r,phantom,"\subset"]
        & \cY \arrow[r,"\Phi_{mL}=\widebar{\varphi}_L"]
            & \cY^{\textup{can}}=\bP(\cC_L)
                & \Phi_{mL}(\cB') \arrow[l,phantom, "\supset"]
\end{tikzcd}
\]

\begin{lem}
\label{l.singular-Y}
    The singular locus of $\cY$ is $S'$ and the germ of $\cY$ at any point $y\in S'$ is analytically isomorphic to 
    \[
    \bC^3/\bZ_2\times \bC^2,
    \]
    where $\bZ_2$ is the cyclic group generated by the action 
    \[
    (x_1,x_2,x_3)\mapsto (-x_1,-x_2,-x_3).
    \]
\end{lem}

\begin{proof}
    According to the construction in \S\,\ref{ss.weak-Fano-model}, the singular locus of $\cY$ is contained in $S'=\nu_1(\widetilde{E}_1)$ as $\cY_1$ is smooth and $\widetilde{E}_1\rightarrow S'$ is a $\bP^2$-bundle such that the restriction of $\sO_{\cY_1}(\widetilde{E}_1)$ to any fiber is isomorphic to $\sO_{\bP^2}(-2)$ by \cite[Proposition 3.15]{HoeringPeternell2024}.  

    { Let $y\in Z_y\subset \cY$ be a general complete intersection of two hypersurfaces passing through $y$ and $Z'_y\coloneqq \nu_1^{-1}(Z_y)$. Then $y$ is an isolated singular point of $Z_y$ and the induced morphism $\nu_1\colon Z'_y\rightarrow Z_y$ is a resolution of singularities so that 
    \[
    E_{y}\coloneqq \nu_1^{-1}(y)\cong \bP^2\quad \text{and}\quad E_y|_{E_y}\cong \sO_{\bP^2}(-2).
    \]
    So $y$ is a quadruple non-Gorenstein point on $Z_y$, which is analytically isomorphic to the quotient of $\bC^3$ under the action $(x_1,x_2,x_3)\mapsto (-x_1,-x_2,-x_3)$ (see for instance \cite[Theorem 1.4.3]{IskovskikhProkhorov1999}).}
\end{proof}

\begin{lem}
    The singular locus of $\cL$ is ${\bf 0}_{S'}$ and the germ of $\cL$ at any point $z\in {\bf 0}_{S'}$ is analytically isomorphic to 
    \[
    \bC^4/\bZ_2\times \bC^2,
    \]
    where ${\bf 0}_{S'}={\bf 0}_{\cL}\cap p^{-1}(S')$ and $\bZ_2$ is the cyclic group generated by the action 
    \[
    (x_1,x_2,x_3,x_4)\mapsto (-x_1,-x_2,-x_3,-x_4).
    \]
    In particular, the variety $\cL$ has only $\bQ$-factorial terminal singularities.
\end{lem}

\begin{proof}
    As $\cY$ is smooth outside $S'$, the divisor $L$ is Cartier outside $S'$. So the singular locus of $\cL$ is contained in $p^{-1}(S')$. We note also that any point $y\in S'$ has index two and $-K_{\cY}$ is the generator of the local divisor class group around $y$ such that $-K_{\cY}\sim K_{\cY}$ around $y$. In particular, as $-K_{\cY}=3L$, we have $L\sim -K_{\cY}\sim K_{\cY}$ around $y$. 

    { By Lemma \ref{l.singular-Y}, the germ of $\cY$ at $y$ is analytically isomorphic to the quotient $\bC^5$ by the involution $i\colon \bC^5\rightarrow \bC^5$ defined as
    \[
    (x_1,x_2,x_3,y_1,y_2) \longmapsto (-x_1,-x_2,-x_3,y_1,y_2).
    \]
    Let $\omega=dx_1\wedge\dots\wedge dx_3\wedge dy_1\wedge dy_2$ be a local generator of the canonical divisor $K_{\bC^5}$ at the origin. Then the involution $i$ induces a natural action $\omega \mapsto -\omega$. Hence the total space of $L\sim K_{\cY}$ around $y$ is isomorphic to the quotient of the canonical bundle
    \[
    K_{\bC^5} \cong \bC^5 \oplus \bC\omega\ni (x_1,x_2,x_3,y_1,y_2,x_4)
    \]
    by the natural $i$-action; that is, the the germ of $\cL$ over $y$ is analytically isomorphic to 
    \[
    K_{\bC^5}/\langle i\rangle \cong \bC^4/\bZ_2\times \bC^2,
    \]
    and its singular locus is just $\{(0,0,0,0)\}\times \bC^2$; that is, it is contained in the zero section ${\bf 0}_{\cL}=\{x_4=0\}$.}
\end{proof}

\begin{prop}
\label{p.QFT-ZV5}
    The symplectic variety $\cZ_X$ admits a unique $\bQ$-factorial terminalization $\mu\colon \cZ'_X\rightarrow \cZ_X$. Set $Y\coloneqq \mu^{-1}(0)$. Then the following statements hold:
    \begin{enumerate}
        \item\label{i1.PV5} The birational map $\cL\rightarrow \cZ_X$ induces an isomorphism $\cL\setminus {\bf 0}_{\cL}\cong \cZ'_X\setminus Y$.
        
        \item\label{i2.PV5} The fiber $Y$ is irreducible and $\mu$ is a symplectic resolution outside $0$.

        \item\label{i3.PV5} $\cZ_{X,\sing}=Z$ and $\cZ_X$ has $A_2$-singularities along the general points of $Z$.

    \end{enumerate}
\end{prop}

\begin{proof}
    Set $E=\textbf{0}_{\cL}$. Then $E\rightarrow \cY$ is an isomorphism and $\sO_{E}(E)\cong \sO_{\cY}(-L)$. On the other hand, as $\varphi_{L}^* K_{\cZ_X} = K_{\cL}-2E$, running a relative $K_{\cL}$-MMP over $\cZ_X$ yields a birational contraction $g\colon \cL\dashrightarrow \cZ'_X$ to a normal variety with $\bQ$-factorial terminal singularities such that it satisfies the following commutative diagram
    \[
    \begin{tikzcd}[column sep=large, row sep=large]
        p^{-1}(\cB') \arrow[r,phantom,"\subset"]
            & \cL \arrow[r,dashed,"g"] \arrow[dr,"\varphi_L" swap]
                & \cZ'_X \arrow[d,"\mu"] 
                    &  \cB'_a \arrow[d] \arrow[l,phantom,"\supset"] \\
            &
                & \cZ_X
                    & Z \arrow[l,phantom,"\supset"]
    \end{tikzcd}
    \]
    where $\cB'_a=g(p^{-1}(\cB'))$. Note that $-E$ is $\varphi_L$-big and $\varphi_L$-nef. In particular, as $K_{\cZ'_X}=g_*{K_{\cL}}$ is $\mu$-nef, the Negativity Lemma implies that $g_*E=0$ and thus $\Exc(g)=E$ as $\Exc(g)\subset E$. Moreover, since $\cL$ is smooth outside $E$ and $Y=g(E)$, it follows that $Y$ is irreducible and $\cZ'_X$ is smooth outside $Y$. This proves \ref{i1.PV5} and \ref{i2.PV5}.
    
    For \ref{i3.PV5}, since $\rho(\cZ'_X)=1$ and $\cB_{a}'$ is contracted by $\mu$, the relative movable cone of $\mu$ is $\bR_{\geq 0}[\cB_a']$. Hence $\cZ'_X$ is the unique $\bQ$-factorial terminalization of $\cZ_X$. Note that $\mu$ is a crepant map, which contracts $\cB'_a$ to $Z$ (cf. Lemma \ref{l.exceptional-locus}). In particular, since $\cZ'_X$ is smooth outside $Y$, the variety $\cZ_X$ is exactly singular along $Z$ $(\ni 0)$. Then the statement follows from Proposition \ref{p.fibres-C'}.
\end{proof}

\begin{prop}
\label{p.SympResV5}
    The variety $\cZ'_X$ is not smooth. In particular, the symplectic orbifold cone $\cZ_X$ does not admit any symplectic resolution.
\end{prop}

\begin{proof}
    Assume that $\cZ'_X$ is smooth. Then $\mu\colon \cZ'_X\rightarrow \cZ_X$ is a symplectic resolution. In particular, the central fiber $Y$ is smooth and $\cZ'_X\cong T^*Y$ by Corollary \ref{c.symp-resolution} and hence $Y$ is a rational homogeneous space as $\dim(\cZ_X)=6$. In particular, we have
    \[
    3=\dim H^0(X,T_X) = \dim S(X)_1 = \dim S(Y)_1 = \dim H^0(Y,T_Y) > 3,
    \]
    which is a contradiction. Finally, we conclude as a symplectic resolution is also a $\bQ$-factorial terminalization. 
\end{proof}

\subsection*{Acknowledgments}
The authors thank Michel Brion and Boris Pasquier for helpful discussions on horospherical varieties. 
{  We are very grateful to the two anonymous referees for their detailed and pertinent reports which help us to improve this paper.}
J.~Liu is supported by the National Key Research and Development Program of China (No. 2021YFA1002300) and the Youth Innovation Promotion Association CAS. Both authors are supported by the CAS Project for Young Scientists in Basic Research (No. YSBR-033) and the NSFC grant (No. 12288201).
	\bibliographystyle{alpha}
	\bibliography{SACS}
\end{document}